\documentclass[final]{siamltex}

\usepackage{epsfig,amssymb,latexsym}
\usepackage{amsfonts,psfrag,amsmath,bbm,color}
\usepackage{cancel}
\usepackage{siunitx}
\usepackage{comment}
\usepackage{mathrsfs}
\usepackage{graphicx}
\usepackage{textcomp}
\usepackage{multirow}
\usepackage{enumerate}
\usepackage{cancel}
\usepackage{algpseudocode}
\usepackage{caption}  
\usepackage{subcaption}
\usepackage{url}
\usepackage{rotating}
\usepackage{bigints}
\usepackage{hyperref}
\usepackage{xcolor}
\usepackage{ulem}
\usepackage{mathrsfs}
\usepackage{placeins}

\newcommand{\bE}{{\bf E}}

\def\grad{{\nabla}}
\usepackage[ruled,vlined]{algorithm2e}
\usepackage{geometry}
\vbadness=\maxdimen

\newcommand{\footremember}[2]{%
    \footnote{#2}
    \newcounter{#1}
    \setcounter{#1}{\value{footnote}}%
}
\newcommand{\footrecall}[1]{%
    \footnotemark[\value{#1}]%
} 
\usepackage{comment}
\usepackage{tikz}
\usetikzlibrary{shapes.geometric, arrows}
\tikzstyle{startstop} = [rectangle, rounded corners, minimum width=1cm, minimum height=1cm,text centered, draw=black]
\tikzstyle{io} = [trapezium, trapezium left angle=70, trapezium right angle=110, minimum width=1cm, minimum height=1cm, text centered, draw=black, fill=blue!30]
\tikzstyle{method} = [rectangle, rounded corners, minimum width=1cm, minimum height =1cm, text centered, draw=black]
\tikzstyle{process} = [rectangle, minimum width=1cm, minimum height=1cm, text centered, draw=black]
\tikzstyle{decision} = [diamond, minimum width=0.5cm, minimum height=0.5cm, text centered, draw=black, fill=green!30]
\tikzstyle{arrow} = [thick,->,>=stealth]

\usepackage{amscd}

\graphicspath{{./figs/}}

\newcommand{\lp}{\left(}
\newcommand{\rp}{\right)}

\newtheorem{remark}{Remark}[section]

\def\PP{{{\rm l}\kern - .15em {\rm P} }}
\def\PN2{{\PP_{N}-\PP_{N-2}}}

\newcommand{\cD}{\mathcal{D}}

\newcommand{\bfeta}{\boldsymbol{\eta}}

\newcommand{\bphi}{\boldsymbol{\varphi}}

\newcommand{\tp}{\tilde{p}}
\newcommand{\tl}{\tilde{l}}

\newcommand{\bif}{\textbf{\textit{f}}\hspace{-0.5mm}}
\newcommand{\tbif}{\tilde{\textbf{\textit{f}}}}
\newcommand{\ba}{\boldsymbol{a}}

\newcommand{\bchi}{\pmb{\chi}}

\newcommand{\be}{\boldsymbol{e}}

\newcommand{\bg}{\textbf{\textit{g}}}

\newcommand{\bH}{\boldsymbol{H}}

\newcommand{\bsU}{\boldsymbol{\mathcal{U}}}

\newcommand{\bL}{\boldsymbol{L}}

\newcommand{\bu}{\boldsymbol{u}}

\newcommand{\bv}{\boldsymbol{v}}
\newcommand{\talpha}{\tilde{\alpha}}

\newcommand{\bV}{\boldsymbol{V}}

\newcommand{\bw}{\boldsymbol{w}}

\newcommand{\btu}{\tilde{\boldsymbol{u}}}

\newcommand{\bx}{\boldsymbol{x}}
\newcommand{\bX}{\boldsymbol{X}}

\newcommand{\by}{\boldsymbol{y}}

\newcommand{\bGamma}{\boldsymbol{\Gamma}}

\newcommand{\deleted}[1]{{}}

\begin{document}
\title{Efficient and Optimally Accurate Numerical Algorithms for Stochastic Turbulent Flow Problems}

\author{
Brandiece N. Berry\footremember{uabm}{D\MakeLowercase{epartment of} M\MakeLowercase{athematics}, U\MakeLowercase{niversity of} A\MakeLowercase{labama at} B\MakeLowercase{irmingham}, AL 35294, USA}%
\and Md Mahmudul Islam\hspace{1mm}\footrecall{uabm}
\and Muhammad Mohebujjaman\footnote{S\MakeLowercase{upported by the} N\MakeLowercase{ational} S\MakeLowercase{cience} F\MakeLowercase{oundation grant} DMS-2425308; C\MakeLowercase{orrespondence: mmohebuj@uab.edu}}\hspace{1mm}\footrecall{uabm}
\and
Neethu Suma Raveendran\footrecall{uabm}
 }

\maketitle

\begin{abstract} 
In this paper, we first propose a filter-based continuous Ensemble Eddy Viscosity (EEV) model for stochastic turbulent flow problems. We then propose a generic algorithm for a family of fully discrete, grad-div regularized, efficient ensemble parameterized schemes for this model. The linearized Implicit-Explicit (IMEX) EEV generic algorithm shares a common coefficient matrix for each realization per time-step, but with different right-hand-side vectors, which reduces the computational cost and memory requirements to the order of solving deterministic flow problems. Two family members of the proposed time-stepping algorithm are analyzed and proven to be stable. It is found that one is first-order and the other is second-order accurate in time for any stable finite element pairs. Avoiding the discrete inverse inequality, the optimal convergence of both schemes is proven rigorously for both 2D and 3D problems.  For appropriately large grad-div parameters, both schemes are unconditionally stable and allow weakly divergence-free elements. Several numerical tests are given for high expected Reynolds number ($\bE[Re]$) problems. The convergence rates are verified using manufactured solutions with $\bE[Re]=10^{3},10^{4},\;\text{and}\; 10^{5}$. For various high  $\bE[Re]$, the schemes are implemented on benchmark problems which includes: A  2D channel flow over a unit step problem, a non-intrusive Stochastic Collocation Method (SCM) is used to examine the performance of the schemes on a 2D  Regularized Lid Driven Cavity (RLDC) problem, and a 3D RLDC problem, and found them perform well.

\end{abstract}



\pagestyle{myheadings}
\thispagestyle{plain}

\markboth{\MakeUppercase{Efficient Alogirthms for Stochastic Turbulent Flow Problems}}{\MakeUppercase{  B. N. Berry, M. M. Islam, M. Mohebujjaman, and N. S. Raveendran}}

\section{Introduction} Laminar or low Reynolds numbers ($Re$) flows are simple and easy to simulate. Complex natural flows are turbulent and ubiquitous.  They are important for many practical applications, e.g., scour formation  around a submerged bridge pier \cite{graf2002flow}, wingtip vortices \cite{clancy1975aerodynamics}, flow around the hull of a submarine \cite{zhang2013simulation} and naval ships \cite{choi2020estimation}, weather prediction \cite{rafalimanana2022optimal}, and flow in the Earth's core \cite{jones2015thermal, olson2013experimental}. Turbulent flows are chaotic in behavior \cite{berselli2006mathematics} and highly ill-conditioned. Accurate simulation of such a flow is one of the most challenging problems due to the limited resolution of discretizations. Even if we could accurately simulate the chaotic behavior of turbulent flow, a significant difficulty in accurate prediction remains in determining the problem data (initial and boundary conditions, viscosity, forcing, etc.) precisely. Indeed, many forecasts that affect our day-to-day existence are derived from dubious turbulent flow simulations. Since the problem data involve random noise, it leads to model turbulent flows using a parameterized Stochastic Navier-Stokes Equations (SNSE), and consider the solution of Navier-Stokes Equations (NSE) at high $Re$ as a realization \cite{barenblatt1996scaling}. In such cases, the predictability of uncertainty quantification problems is likely to increase.

Let $\mathcal{D}\subset \mathbb{R}^d\ (d\in\{2,3\})$ be a convex polygonal or polyhedral physical domain with boundary $\partial\mathcal{D}$. A complete probability space is denoted by $(\Omega,\mathcal{F},P)$ with $\Omega$ the set of outcomes, $\mathcal{F}\subset 2^\Omega$ the $\sigma$-algebra of events, and $P:\mathcal{F}\rightarrow [0,1]$ represents a probability measure. We consider the time-dependent, dimensionless, viscoresistive, incompressible SNSE turbulent flow problems for homogeneous Newtonian fluids, which are governed by the following non-linear stochastic Partial Differential Equations (PDE) \cite{gunzburger2019evolve}:\setlength{\abovedisplayskip}{0pt} \setlength{\abovedisplayshortskip}{0pt}
\setlength{\belowdisplayskip}{0pt} \setlength{\belowdisplayshortskip}{0pt}\vspace{-1ex}\begin{align}
    \btu_{t}+\btu\cdot\nabla \btu-\nabla\cdot\left(\nu(\bx,\omega)\nabla\btu\right)+\nabla \tp &=  \tbif(\bx,t,\omega), \hspace{2mm}\text{in}\hspace{2mm}\mathcal{D} \times (0,T]\times\Omega,\label{momentum}\\
\nabla\cdot \btu & = 0, \hspace{14.5mm}\text{in}\hspace{2mm}\mathcal{D} \times (0,T]\times\Omega, \\
\btu(\bx,t,\omega)&=\tilde{\bg}(\bx,\omega),\hspace{5mm}\text{in}\hspace{2.5mm} \partial\cD\times(0,T]\times\Omega,\label{bc-condition}\\
\btu(\bx,0,\omega)& = \btu^0(\bx,\omega),\hspace{3.7mm}\text{in}\hspace{2mm}\mathcal{D}\times\Omega.\label{nse-initial}
\end{align}
The simulation end time is represented by $T>0$, $\bx$ the spatial variable, and $t$ the time variable. The random viscosity $\nu(\bx,\omega)$, external forcing function $\tbif(\bx,t,\omega)$, prescribed Dirichlet boundary condition $\tilde{\bg}(\bx,\omega)$, and initial condition $\btu^0(\bx,\omega)$ are modeled as a random field with $\omega\in\Omega$ for a sample space $\Omega$. The unknown quantities are the velocity field $\btu(\bx,t,\omega)\in\mathbb{R}^d$, and the modified pressure $\tp(\bx,t,\omega)\in\mathbb{R}$ with mean zero. 


 Numerical models for \eqref{momentum}-\eqref{nse-initial} encounter several challenges. For Uncertainty Quantification (UQ), it is common practice to draw repeated random samples of the model inputs (in the Monte Carlo Method (MCM) or SCM), run the complete simulations for each realization, and compute the average or weighted average Quantity of Interest (QoI), e.g., energy, lift, and drag. This leads the above problem \eqref{momentum}-\eqref{nse-initial} into computing the flow ensemble solving the following system: 
\begin{align}
\btu_{j,t}+\btu_j\cdot\nabla \btu_j-\nabla\cdot\left(\nu_j(\bx) \nabla \btu_j\right)+\nabla \tp_j &=  \tbif_j(\bx,t), \hspace{2mm}\text{in}\hspace{2mm}\cD \times (0,T], \label{gov1}\\
\nabla\cdot \btu_j & = 0, \hspace{11.7mm}\text{in}\hspace{2.mm}\cD \times (0,T],\label{gov2}\\
\btu_j(\bx,t)&=\tilde{\bg}_j(\bx),\hspace{5mm}\text{in}\hspace{2mm} \partial\cD\times(0,T],\label{gov3}\\
\btu_j(\bx,0)& =\btu_j^0(\bx),\hspace{5mm}\text{in}\hspace{2mm}\mathcal{D},\label{gov4}
\end{align}
where $\btu_j$, and $\tp_j$, denote the velocity, and pressure solutions, respectively, for each $j=1,2,\cdots,J$, corresponding to distinct kinematic viscosity $\nu_j$, and/or body forces $\tbif_j$,  and/or the initial conditions $\btu_j^0$, and/or boundary conditions $\tilde{\bg}_j$. Here $J$ the total number of realization, typically large. To make the analysis simple, set $\tilde{\bg}_j=\textbf{0}$. 

For each realization, to have a resolved scale flow, at each time-step, a Direct Numerical Simulation (DNS), is required to be done with an extremely large number of degrees of freedom (dof), e.g., for flow around a moderate size airplane dof $\approx 2.7\times 10^{16}$. The required number of mesh points is $O(Re^{\frac{9}{4}})$ in 3D. This is due to the presence of a wide range of length scales and a hierarchy of eddies, which is known from the Kolmogorov theory \cite{kolmogorov1941local, kolmogorov1941degeneration}. To handle this issue, We follow the standard Large Eddy Simulation (LES) approach \cite{layton2012approximate,rebollo2014mathematical}, and assume the filtering operator commutes with differentiation. To this end, we introduce $\btu_j=\bu_j+\bw_j$, and $\tp_j=p_j+r_j$ in \eqref{gov1}-\eqref{gov2}, and finally, take filter operator on both sides, then we have the following large scale system \cite{berselli2006mathematics,durbin2011statistical,wilcox2006turbulence} 
\setlength{\abovedisplayskip}{0pt} \setlength{\abovedisplayshortskip}{0pt}
\setlength{\belowdisplayskip}{0pt} \setlength{\belowdisplayshortskip}{0pt}\begin{align}
	\bu_{j,t}+\bu_j\cdot\nabla \bu_j+\nabla\cdot(\overline{\bw_j\bw_j^T})-2\nabla\cdot\left(\nu_j(\bx) \nabla \bu_j\right)+\nabla p_j &=  \bif_j, \label{gov1-rans}\\
	\nabla\cdot \bu_j & = 0,\label{gov2-rans}
\end{align}where  $\bu_j:=\overline{\btu}_j$, $p_j:=\overline{\tp}_j$, and $\bif_j:=\overline{\tbif}_j$ for which $\overline{\bw}_j=\textbf{0}$, and $\overline{r}_j=0$. The new ensemble system \eqref{gov1-rans}-\eqref{gov2-rans} is open;  thus, the Reynolds Stress Tensor (RST) \cite{durbin2011statistical} $\overline{\bw_j\bw_j^T}$  is estimated and forms a closed-form systems. It is popular to model RST as eddy viscosity, e.g., $k$-$\epsilon$ \cite{mohammadi1993analysis}, LES \cite{berselli2006mathematics}, Smagorinsky \cite{jiang2016algorithms,siddiqua2023numerical,smagorinsky1963general} models, and EEV \cite{jiang2015higher, jiang2015numerical,mohebujjaman2024efficient,mohebujjaman2017efficient,mohebujjaman2022efficient,raveendran2025efficient}, or using the Boussinesq assumption \cite{boussinesq1877essai} to have closed models. Thus, we approximate $\nabla\cdot(\overline{\bw_j\bw_j^T})\approx-\nabla\cdot\left(2\nu_T\nabla\bu_j\right)$ \cite{berselli2006mathematics,durbin2011statistical,layton2012approximate}. Here, the turbulent EEV coefficient $\nu_T$ is given by the Kolmogorov-Prandtl relation \cite{jiang2015numerical} $\nu_T:=\mu l^2\delta$ where $\mu$ is a tuning parameter, the Prandtl mixing length\vspace{-2mm} \begin{align}
l:=\sqrt{\sum_{j=1}^J|\bu_j^{'}|^2} \label{Prandtl- mixing-length}
\end{align}is a scalar quantity,  $|\cdot|$ denotes length of a vector, the velocity fluctuation $\bu_j^{'}:=\bu_j-\frac{1}{J}\sum\limits_{j=1}^J\bu_j$, the radius of the filter $\delta=\Delta t$, and $\Delta t$ is the time-step size. This approximation models the effect of fluctuations on large scales. 

The simulation of \eqref{gov1-rans}-\eqref{gov2-rans} requires a huge storage and a long simulation time. Because their standard variational discretization (linearized or not) eventually leads to a saddle point linear system whose coefficient matrix varies with the index $j$ at each time-step as: For $j=1,2,\cdots,J$\begin{align}
\begin{pmatrix}\mathbb{A}_j & \mathbb{B}^T\\\mathbb{B} & \mathcal{O}\end{pmatrix}\bigg(\begin{matrix}
\textbf{U}_j\\
     \textbf{P}_j
\end{matrix}\bigg)=\bigg(\begin{matrix}
    \textbf{F}_j\\\textbf{G}_j
\end{matrix}\bigg),\label{sparse-system}
 \end{align}
 where $\textbf{U}_j$, and $\textbf{P}_j$ are the nodal vectors of velocity, and pressure fields, respectively, and $\textbf{F}_j$ and $\textbf{G}_j$ are for the appropriate right-hand-side vector,
 $\mathbb{A}_j$ is the block-matrix corresponding to $\textbf{U}_j$, $\mathbb{B}$ represents the gradient operator, $\mathbb{B}^T$ its adjoint, and $\mathcal{O}$ represents a zero block-matrix. Therefore, the total cost is equal to $J\times$ One NSE simulation. To reduce this huge simulation time and memory requirement, Jiang and Layton \cite{JL14} proposed a breakthrough efficient scheme which leads to a system at each time-step below: For $j=1,2,\cdots,J$\begin{align}
\begin{pmatrix}\mathbb{A} & \mathbb{B}^T\\\mathbb{B} & \mathcal{O}\end{pmatrix}\bigg(\begin{matrix}
\textbf{U}_j\\
     \textbf{P}_j
\end{matrix}\bigg)=\bigg(\begin{matrix}
    \textbf{F}_j\\\textbf{G}_j
\end{matrix}\bigg),\label{sparse-system2}
 \end{align}
 where the coefficient matrix is independent of the index $j$, moreover, it can take the advantage of the block linear system as
 \begin{align}
\begin{pmatrix}\mathbb{A} & \mathbb{B}^T\\\mathbb{B} & \mathcal{O}\end{pmatrix}\bigg(\begin{matrix}
\textbf{U}_1\\
     \textbf{P}_1
\end{matrix}
\bigg\rvert\begin{matrix}
\textbf{U}_2\\
\textbf{P}_2
\end{matrix}\bigg\rvert\begin{matrix}
    \cdots\\
    \cdots
\end{matrix}\bigg\rvert\begin{matrix}
\textbf{U}_J\\
\textbf{P}_J
\end{matrix}\bigg)=\bigg(\begin{matrix}
    \textbf{F}_1\\\textbf{G}_1
\end{matrix}\bigg\rvert\begin{matrix}
    \textbf{F}_2\\\textbf{G}_2
\end{matrix}\bigg\rvert\begin{matrix}
    \cdots\\
    \cdots
\end{matrix}\bigg\rvert\begin{matrix}
    \textbf{F}_J\\\textbf{G}_J
\end{matrix}\bigg).\label{sparse-system-block}
 \end{align}
 Ensemble algorithms having this efficient feature have been implemented in many areas, e.g., NSE \cite{gunzburger2019efficient,JL14,jiang2024second,raveendran2025efficient,yuan2020ensemble}, Boussinesq Equations \cite{jiang2024ensemble}, and Magnetohydrodynamics (MHD) \cite{jiang2018efficient, Mohebujjaman2022High,MR17}.
 For parameterized NSE flow problems an efficient first and second order time-stepping ensemble algorithm is proposed in \cite{gunzburger2019efficient} and \cite{jiang2024second}, respectively. 

Commonly used Taylor-Hood (TH) \cite{wieners2003taylor} elements are not pointwise divergence-free; thus, the divergence error propagates and reduces the fidelity of the solution. Grad-div stabilization with an appropriately large coefficient can enforce the divergence-free constraint pointwise. Strongly divergence-free elements like Scott-Vogelius (SV) \cite{scott1985norm} require higher computational cost.


 The viscosity $\nu$ is the most critical and sensitive key parameter as it determines the flow characteristics, but cannot be measured perfectly; consequently, the models exhibit the ``butterfly effect" \cite{lorenz1963deterministic} in their predictions. The parameterized numerical models accounts the uncertainty in the viscosity parameter to achieve high-fidelity solutions, in addition to the presence of random noise in other input data. To the best of our knowledge, no parameterized EEV efficient scheme has been proposed yet for NSE problems. An efficient parameterized time-stepping first-order and second-order algorithm is proposed in \cite{gunzburger2019efficient} and \cite{Gunzburger2019Second-Order}, respectively, for laminar flow problems (without modeling RST term in \eqref{gov1-rans}). Recently, several efficient EEV schemes are proposed for non-parameterized NSE \cite{jiang2015higher,jiang2015numerical,jiang2015analysis,jiang2016algorithms,jiang2020foundations} problems, however, avoiding the convergence proofs, the order of convergence are exhibited numerically. They provided numerical experiments for laminar flow problems. Several efficient EEV schemes have been proposed for laminar Magnetohydrodynamic (MHD) flow problems in \cite{mohebujjaman2017efficient, mohebujjaman2022efficient} in recent years, and they have been proven to be sub-optimal convergent for 3D problems.  Penalty-projection based EEV schemes are proposed for laminar NSE in \cite{raveendran2025efficient} and for laminar MHD in \cite{mohebujjaman2024efficient} where the convergence of the discrete schemes to the respective (temporally sub-optimal) coupled schemes are given.  


In this work, we propose a parameterized, efficient, grad-div regularized, fully discrete, time-stepping generic algorithm for the variational formulation of \eqref{gov1-rans}-\eqref{Prandtl- mixing-length}. The algorithm is elegantly designed to possess the property in \eqref{sparse-system2}, which saves a huge amount of computational cost and memory requirement.  Moreover, it is a linearized IMEX algorithm based on the Backward Differentiation Formula (BDF), so that at each time-step, there is no need to solve a system of nonlinear algebraic equations, which dramatically reduces the number of arithmetic operations. Depending on the problem, an appropriately large grad-div parameter enforces the discrete divergence-free constraints, allowing the use of a weakly divergence-free element that requires fewer degrees of freedom (dof) compared to a pointwise divergence-free element. We then analyze and test two members of the proposed generic stochastic turbulent flow algorithm: Linearized Backward-Euler (BE) and linearized BDF-2 schemes.  We carefully analyze and have found that the use of discrete inverse inequality in EEV terms leads to sub-optimal convergence in \cite{mohebujjaman2017efficient, mohebujjaman2022efficient}. In our analysis, by avoiding the discrete inverse inequality, we prove that stable time-stepping schemes are optimally convergent both in space and time. Specifically, the BE-EEV scheme is first-order accurate, while the BDF-2-EEV scheme is second-order accurate in time. Numerical experiments are presented to verify the predicted convergence rates for the expected viscosities $\bE[\nu]=10^{-3}$, $10^{-4}$, and $10^{-5}$ using a manufactured solution. The proposed schemes have been tested on 2D/3D benchmark problems with high $\bE[Re]$ to evaluate their efficiency and performance, which include: A 2D rectangular channel flow past a unit step problem, a 2D RLDC problem in conjunction with non-intrusive modular SCM, and a 3D RLDC problem. The dimension-independent codes for these schemes are written on the Deal.II \cite{dealII93} finite element library.

To the best of our knowledge, the proposed filter-based continuous ensemble model in \eqref{gov1-rans}-\eqref{Prandtl- mixing-length}, the fully discrete EEV generic algorithm, and the analysis of its two members are novel contributions to the study of turbulent SNSE problems.

  The rest of the paper is organized as follows: Notations and mathematical preliminaries are presented in Section \ref{notation-prelims}, which are used throughout the paper. In Section \ref{fully-discrete-scheme}, we present and analyze two fully discrete, linearized, efficient, and grad-div regularized schemes. A series of numerical experiments is given in Section \ref{numerical-experiment}, which support the theory, combine SCM with the proposed schemes, and implement them on several benchmark problems. Finally, conclusions are drawn and future research directions are given in Section \ref{conclusion}.

\section{Notation and preliminaries}\label{notation-prelims}

Let $\cD\subset \mathbb{R}^d\ (d=2,3)$ be a convex polygonal or polyhedral domain in $\mathbb{R}^d(d=2,3)$ with boundary $\partial\cD$. The usual $L^2(\cD)$ norm and inner product are denoted by $\|.\|$ and $(.,.)$, respectively. Similarly, the $L^p(\cD)$ norms and the Sobolev $W_p^k(\cD)$ norms are $\|.\|_{L^p}$ and $\|.\|_{W_p^k}$, respectively for $k\in\mathbb{N},\hspace{1mm}1\le p\le \infty$. Sobolev space $W_2^k(\cD)$ is represented by $H^k(\cD)$ with norm $\|.\|_k$. The vector-valued spaces are $$\bL^p(\cD)=(L^p(\cD))^d, \hspace{1mm}\text{and}\hspace{1mm}\bH^k(\cD)=(H^k(\cD))^d.$$
For $\bX$ being a normed function space in $\cD$, $L^p(0,T;\bX)$ is the space of all functions defined on $(0,T]\times\cD$ for which the following norm\begin{align*} \|\bu\|_{L^p(0,T;\bX)}=\lp\int_0^T\|\bu\|_{\bX}^pdt\rp^\frac{1}{p},\hspace{2mm}p\in[1,\infty)
\end{align*}
is finite. For $p=\infty$, the usual modification is used in the definition of this space. The natural function spaces for our problem are
\begin{align*}
    \bX:&=\bH_0^1(\cD)=\{\bv\in \bL^2(\cD) :\nabla \bv\in L^2(\cD)^{d\times d}, \bv=0 \hspace{2mm} \mbox{on}\hspace{2mm}   \partial \cD\},\\
    Q:&=L_0^2(\cD)=\{ q\in L^2(\cD): \int_\cD q\hspace{1mm}d\bx=0\}.
\end{align*}

Recall the Poincare inequality holds in $\bX$: There exists $C$ depending only on $\cD$ satisfying for all $\bphi\in \bX$,
\[
\| \bphi \| \le C \| \nabla \bphi \|.
\]

We define the skew symmetric trilinear form $b^*:\bX\times \bX\times \bX\rightarrow \mathbb{R}$ by
	\[
	b^*(\bu,\bv,\bw):=\frac12(\bu\cdot\nabla \bv,\bw)-\frac12(\bu\cdot\nabla \bw,\bv). 
	\]
    
 By the divergence theorem \cite{jiang2015higher}, it can be shown \begin{align}
     b^*(\bu,\bv,\bw)=(\bu\cdot\nabla \bv,\bw)+\frac12(\nabla\cdot\bu,\bv\cdot\bw).\label{trilinear-identitiy}
 \end{align}
 
 Recall from \cite{L08, lee2011error, linke2017connection} that for any $\bu,\bv,\bw\in 
		\bX$
	\begin{align}
		b^*(\bu,\bv,\bw)&\leq C(\cD)\|\nabla \bu\|\|\nabla \bv\|\|\nabla \bw\|,\label{nonlinearbound}
	\end{align}	
 and additionally, if $\bv\in \bL^\infty(\cD)$, and $\nabla\bv\in\bL^3(\cD)$, then 
 \begin{align}
    b^*(\bu,\bv,\bw)\leq C(\cD)\|\bu\|\left(\|\nabla\bv\|_{L^3}+\|\bv\|_{L^\infty}\right)\|\nabla\bw\|. \label{nonlinearbound3}
 \end{align}
 
 The following basic inequality will be used
\begin{align}
\|\bu\cdot\nabla\bv\|&\le\||\bu|\nabla\bv\|.\label{basic-ineq}
\end{align}

The space of divergence free functions is given by
$$\bV:=\{\bv\in\bX:(\nabla\cdot\bv,q)=0,\forall q\in Q\}.$$

The conforming finite element spaces are denoted by $\bX_h\subset \bX$ and  $Q_h\subset Q$, and we assume a regular triangulation $\tau_h(\cD)$, where $h$ is the maximum triangle diameter.   We assume that $(\bX_h,Q_h)$ satisfies the usual discrete inf-sup condition
\begin{eqnarray}
\inf_{q_h\in Q_h}\sup_{\bv_h\in \bX_h}\frac{(q_h,\grad\cdot \bv_h)}{\|q_h\|\|\grad \bv_h\|}\geq\beta>0,\label{infsup}
\end{eqnarray}
where $\beta$ is independent of $h$. 

The discretely divergence free subspace of $\bX_h$ is $\bV_h$ which is conforming to $\bV$, i.e., $\bV_h\subset\bV$, where$$\bV_h:=\{\bv_h\in\bX_h:\left(\nabla\cdot\bv_h,q_h\right)=0,\forall q_h\in Q_h\}.$$

We have the following approximation properties in $(\bV_h,Q_h)$ for piecewise polynomials of degree $(k,k-1)$ \cite{mohebujjaman2022efficient}\begin{align}
\| \bu- P^{L^2}_{\bV_h}(\bu)   \|&\leq Ch^{k+1}|\bu|_{k+1},\hspace{2mm}\bu\in \bH^{k+1}(\cD),\label{App}\\
 \| \nabla (\bu- P^{L^2}_{\bV_h}(\bu)  ) \|&\leq Ch^{k}|\bu|_{k+1},\hspace{6mm}\bu\in \bH^{k+1}(\cD),\label{AppPro3}
\end{align}
where $|\cdot|_r$ denotes the $\bH^r$ seminorm and $P^{L^2}_{\bV_h}(\bu)$ is the $L^2$ projection of $\bu$ into $\bV_h$. The following lemma for the discrete Gr\"onwall inequality was given in \cite{HR90}.
\begin{lemma}\label{dgl}
		Let $\Delta t$, $\mathcal{E}$, $a_n$, $b_n$, $c_n$, $d_n$ be non-negative numbers for $n=1,\cdots, M$ such that
		$$a_M+\Delta t \sum_{n=1}^Mb_n\leq \Delta t\sum_{n=1}^{M-1}{d_na_n}+\Delta 
		t\sum_{n=1}^Mc_n+\mathcal{E}\hspace{3mm}\mbox{for}\hspace{2mm}M\in\mathbb{N},$$
		then for all $\Delta t> 0,$
		$$a_M+\Delta t\sum_{n=1}^Mb_n\leq \mbox{exp}\left(\Delta t\sum_{n=1}^{M-1}d_n\right)\lp\Delta 
		t\sum_{n=1}^Mc_n+\mathcal{E}\rp\hspace{2mm}\mbox{for}\hspace{2mm}M\in\mathbb{N}.$$
	\end{lemma}
We assume $\nu_j(\bx)\in L^\infty(\cD)$, and $\nu_j(\bx)\ge\nu_{j,\min}>0$, where $\nu_{j,\min}=\min\limits_{\bx\in\cD}\nu_j(\bx)$, for $j=1,2,\cdots,J$. Also, define $\Bar{\nu}_{\min}:=\min\limits_{\bx\in\cD}\Bar{\nu}(\bx)$.

 \section{Family of efficient fully discrete  EEV algorithms for turbulent SNSE}\label{fully-discrete-scheme}
In this section, we propose a generic fully discrete, parameterized, efficient, grad-div regularized, linear extrapolated IMEX finite element time-stepping algorithm for the variational formulation of \eqref{gov1-rans}-\eqref{Prandtl- mixing-length}.  

We consider a uniform time-step size $\Delta t$ and let $t_n=n\Delta t$ for $n=0, 1, \cdots$. Assume that the initial condition(s) are in $\bH^2(\cD)\cap\bV$, and present the generic extrapolated BDF EEV scheme in Algorithm \ref{coupled-alg-com}. \\
\begin{algorithm}[H]\label{coupled-alg-com}
  \caption{Family of fully discrete linearized BDF EEV scheme} Given time-step $\Delta t>0$, end time $T>0$, set $M=T/\Delta t$ and if $\bif_{j}\in$ $ L^\infty\left( 0,T;\bH^{-1}(\cD)\right)$ for $j=1,2,\cdots\hspace{-0.35mm},J$. Compute: Find $(\bu_{j,h}^{n+1}, p_{j,h}^{n+1})\in \bX_h\times Q_h$ satisfying, for all $(\bchi_{h},q_{h})\in \bX_h\times Q_h$:
 \begin{align}
\frac{\beta}{\Delta t}&\Big(\bu_{j,h}^{n+1},\bchi_{h}\Big)+b^*\big(\hspace{-1mm}<\bu_h>^n, \bu_{j,h}^{n+1},\bchi_{h}\big)+\big(\Bar{\nu}\nabla \bu_{j,h}^{n+1},\nabla\bchi_{h}\big)+(\gamma\nabla\cdot\bu_{j,h}^{n+1}-p_{j,h}^{n+1},\nabla\cdot\bchi_{h})\nonumber\\&+\left( 2\nu_T^h\nabla \bu_{j,h}^{n+1},\nabla\bchi_h\right)= \big(\tilde{\bif}_{j},\bchi_{h}\big)-b^*(\bu_{j,h}^{'n}, \bsU_{j,h}^n,\bchi_{h})-\big(\nu_j^{'}\nabla \bsU_{j,h}^n,\nabla\bchi_{h}\big),
\label{couple-eqn-1-new}\\\nonumber\\&\big(\nabla\cdot\bu_{j,h}^{n+1},q_{h}\big)=0,\label{couple-incompressibility-new}
\end{align}
\end{algorithm}
\noindent where $\bu_{j,h}^n$, and $p_{j,h}^{n}$ denote approximations of $\bu_j(\cdot,t_n)$, and $p_j(\cdot,t_n)$, respectively. The ensemble mean 
$$<\bu_h>^n:=\sum_{j=1}^J\bsU_{j,h}^n,\text{ and }\;\Bar{\nu}:=\frac{1}{J}\sum\limits_{j=1}^{J}\nu_j(\bx),$$ 
where $\bsU_{j,h}^n$ represents an extrapolated solenoidal velocity. The viscosity fluctuation is defined as $\nu_j^{'}:=\nu_j-\Bar{\nu}$. The discrete EEV coefficient \cite{jiang2015higher,jiang2015analysis,jiang2015numerical,Mohebujjaman2022High, mohebujjaman2024efficient, MR17,mohebujjaman2022efficient} is defined as \begin{align} \nu_T^h:=\mu\Delta t(l_h^n)^2,\hspace{0.5mm}\text{ where}\hspace{1mm}l_h^n:=\sqrt{\sum_{j=1}^J|\bu_{j,h}^{'n}|^2}. \label{eddy-viscosity-new}
\end{align} 

The grad-div stabilization parameter $\gamma>0$ is used to enforce the discrete divergence-free constraint pointwise \cite{jenkins2014parameter,linke2011convergence} for the TH elements. It allows for the avoidance of SV the elements. The SV elements require bary-centered refined meshes, which are only defined for triangular/tetrahedrahe meshes (that is, does not define for quadrilateral/hexahedral meshes) and come with higher dof. The EEV term provides long-time stability for unresolved meshes in convection-dominated flows. The algorithm is efficient because it possesses the feature given in \eqref{sparse-system-block}.  Therefore, it allows for saving a huge computational time and computer memory. For a robust and high fidelity solution, computation of ensemble average solution is popular in many applications such as NSE \cite{pei2024variable}, surface data assimilation \cite{fujita2007surface}, MHD \cite{jiang2018efficient}, porous media flow \cite{jiang2021artificial}, weather forecasting \cite{L05,LP08}, spectral methods \cite{LK10}, sensitivity analyses \cite{MX06}, and hydrology \cite{GG11}.


 Depending on the choice of BDF family implicit schemes and appropriate extrapolation for linearization, the proposed generic algorithm provides both single and linear multi-step methods. 
 For example, we can derive efficient linearized EEV schemes as below
 
 \begin{itemize}
    \item BDF-1 or BE: $\beta:=1,\; \bsU_{j,h}^n:=\bu_{j,h}^n,\; \bu_{j,h}^{'n}:=\bu_{j,h}^n-<\bu_h>^n
    $ and $\tilde{\bif}_{j}:=\bif_{j}(t_{n+1})+\frac{1}{\Delta t}\bu_{j,h}^n.$\\
    \item BDF-2: $\beta:=\frac{3}{2},\; \bsU_{j,h}^n:=2\bu_{j,h}^n-\bu_{j,h}^{n-1},\; \bu_{j,h}^{'n}:=2\bu_{j,h}^n-\bu_{j,h}^{n-1}-<\bu_h>^n$, and $\tilde{\bif}_{j}:=\bif_{j}(t_{n+1})+\frac{2}{\Delta t}\bu_{j,h}^n-\frac{1}{2\Delta t}\bu_{j,h}^{n-1}.$\\
    \item BDF-3: $\beta:=\frac{11}{6},\; \bsU_{j,h}^n:=3\bu_{j,h}^n-3\bu_{j,h}^{n-1}+\bu_{j,h}^{n-2},\; \bu_{j,h}^{'n}:=3\bu_{j,h}^n-3\bu_{j,h}^{n-1}+\bu_{j,h}^{n-2}-<\bu_h>^n$, and $\tilde{\bif}_{j}:=\bif_{j}(t_{n+1})+\frac{3}{\Delta t}\bu_{j,h}^n-\frac{3}{2\Delta t}\bu_{j,h}^{n-1}+\frac{1}{3\Delta t}\bu_{j,h}^{n-2}.$\\
    \item And so on.
\end{itemize} 

The optimal temporal accuracy for BDF-1 is $O(\Delta t)$, for BDF-2 is $O(\Delta t^2)$, for BDF-3 is $O(\Delta t^3)$, and so on. The lower-order methods are less accurate but require fewer arithmetic operations, and are thus popular. On the other hand,  the higher the order, the higher the accuracy of the solutions, but it requires a higher number of arithmetic operations and suffers from the stability issue and involvement of the round-off errors.

In this paper, we analyze and test two members of the above BDF family: the BE-EEV and BDF-2-EEV.
\begin{lemma}\label{lemma1} For $j=1,2,\cdots, J$,  there exists a constant $K^*>0$ such that $\|\bu_{j,h}^n\|_{L^\infty}\le K^*$, $\forall n\in\mathbb{N}$, where $\bu_{j,h}^n$ is a solution of the Algorithm \ref{coupled-alg-com}.
\end{lemma}
\begin{proof}
    The straightforward proof is similar as in the Lemma 3.5 given in \cite{mohebujjaman2024decoupled}.
\end{proof}

 \subsection{Efficient linearized  BE-EEV scheme} To have an efficient linearized BE-EEV scheme from Algorithm \ref{coupled-alg-com}, we consider the data stated above for BDF-1. That is, $$\beta=1,\; \bsU_{j,h}^n=\bu_{j,h}^n,\; \bu_{j,h}^{'n}=\bu_{j,h}^n-<\bu_h>^n,\;\text{and}\;\tilde{\bif}_{j}=\bif_{j}(t_{n+1})+\frac{1}{\Delta t}\bu_{j,h}^n.$$
  To simplify the notation, denote $\alpha_j:=\Bar{\nu}_{\min}-\|\nu_j^{'}\|_\infty$, for $j=1,2,\cdots\hspace{-0.35mm}, J$. We assume that the data does not have outlier and observations are close enough to the mean so that $\alpha_j>0$ holds, that is, $\frac{\|\nu_j^{'}\|_\infty}{\Bar{\nu}_{\min}}<1$. We provide the stability and convergence theorems of the BE-EEV scheme below.
\begin{theorem}(Stability of the BE-EEV scheme)\label{stability-BE-EEV}
Suppose $\bif_{\hspace{0.5mm}j}\in L^\infty(0,T;\bH^{-1}(\cD))$, and $\bu_{j,h}^0=\bu_j(\bx,0)\in\bH^2(\cD)$ for all $j=1,2,
\cdots, J$, then the solutions of BE-EEV scheme (Algorithm \ref{coupled-alg-com}) are stable: Given \begin{align}
    \Delta t\le\min_{\substack{1\le j\le J \\ 1\le n\le M}}\frac{C\alpha_j}{\|\nabla\cdot\bu_{j,h}^{'n}\|^2_{L^\infty}},\label{time-step-BE}
\end{align} if $\alpha_j>0$, choose $\mu\ge \frac12$, we then have the following stability bound:
\begin{align}
\|\bu_{j,h}^{M}\|^2+\alpha_j\Delta t\sum_{n=1}^{M}\|\nabla\bu_{j,h}^n\|^2+2 \gamma\Delta t\sum_{n=1}^{M}\|\nabla\cdot\bu_{j,h}^{n}\|^2\nonumber\\\le \|\bu_{j,h}^{0}\|^2+\Bar{\nu}_{\min}\Delta t\|\nabla \bu_{j,h}^{0}\|^2+\frac{2\Delta t}{\alpha_j}\sum_{n=1}^{M}\|\bif_{\hspace{0.5mm}j}(t_{n})\|_{-1}^2.\label{stability-couple-alg}
\end{align}
\end{theorem}

\begin{proof}
	Taking $\bchi_{h}=\bu_{j,h}^{n+1}\in\bX_h$ and $q_{h}=p^{n+1}_{j,h}\in Q_h$ in \eqref{couple-eqn-1-new} and \eqref{couple-incompressibility-new}, respectively, using $b^*\big(\hspace{-1mm}<\bu_h>^n, \bu_{j,h}^{n+1},\bu_{j,h}^{n+1}\big)=0$, we obtain
	\begin{align}
		&\Big(\frac{\bu_{j,h}^{n+1}-\bu_{j,h}^n}{\Delta t},\bu_{j,h}^{n+1}\Big)+\|\Bar{\nu}^{\frac12}\nabla \bu_{j,h}^{n+1}\|^2+\gamma\|\nabla\cdot\bu_{j,h}^{n+1}\|^2+\left(2\nu_T^h\nabla \bu_{j,h}^{n+1},\nabla \bu_{j,h}^{n+1}\right)\nonumber\\&= (\bif_{j}(t_{n+1}),\bu_{j,h}^{n+1})-(\bu_{j,h}^{'n}\cdot\nabla \bu_{j,h}^n,\bu_{j,h}^{n+1})-(\nu_j^{'}\nabla \bu_{j,h}^{n},\nabla\bu_{j,h}^{n+1}).
	\end{align}

	Using polarization identity and $\left(2\nu_T^h\nabla \bu_{j,h}^{n+1},\nabla \bu_{j,h}^{n+1}\right)=2\mu\Delta t\|l_h^n\nabla \bu_{j,h}^{n+1}\|^2$, we have
	\begin{align}
		&\frac{1}{2\Delta t}\Big(\|\bu_{j,h}^{n+1}\|^2-\|\bu_{j,h}^n\|^2+\|\bu_{j,h}^{n+1}-\bu_{j,h}^n\|^2\Big)+\|\Bar{\nu}^\frac12\nabla \bu_{j,h}^{n+1}\|^2+\gamma\|\nabla\cdot\bu_{j,h}^{n+1}\|^2\nonumber\\&+2\mu\Delta t\|l_h^n\nabla \bu_{j,h}^{n+1}\|^2= (\bif_{j}(t_{n+1}),\bu_{j,h}^{n+1})-b^*(\bu_{j,h}^{'n}, \bu_{j,h}^n,\bu_{j,h}^{n+1})-(\nu_j^{'}\nabla \bu_{j,h}^{n},\nabla\bu_{j,h}^{n+1}).\label{before-nl-bound-new}
	\end{align}
    
  Applying Cauchy-Schwarz and Young's inequalities on the forcing term, yields
  \begin{align*}
      (\bif_{j}(t_{n+1}),\bu_{j,h}^{n+1})\le\|\bif_{j}(t_{n+1})\|_{-1}\|\nabla\bu_{j,h}^{n+1}\|\le \frac{\alpha_j}{4}\|\nabla\bu_{j,h}^{n+1}\|^2+\frac{1}{\alpha_j}\|\bif_{j}(t_{n+1})\|_{-1}^2.
  \end{align*}
  
  We rewrite the trilinear form in \eqref{before-nl-bound-new}, use identity \eqref{trilinear-identitiy}, Cauchy-Schwarz, H\"older's, Poincar\'e, \eqref{basic-ineq}, \eqref{eddy-viscosity-new}, and Young's inequalities, to have
\begin{align*}
b^*(\bu_{j,h}^{'n}, &\bu_{j,h}^n,\bu_{j,h}^{n+1})=b^*(\bu_{j,h}^{'n},\bu_{j,h}^{n+1} ,\bu_{j,h}^{n+1}-\bu_{j,h}^n)\nonumber\\&=(\bu_{j,h}^{'n}\cdot\nabla\bu_{j,h}^{n+1} ,\bu_{j,h}^{n+1}-\bu_{j,h}^n)+\frac12\left(\nabla\cdot\bu_{j,h}^{'n},\bu_{j,h}^{n+1}\cdot(\bu_{j,h}^{n+1}-\bu_{j,h}^n)\right)\nonumber\\
   &\le\|\bu_{j,h}^{'n}\cdot\nabla
\bu_{j,h}^{n+1}\|\|\bu_{j,h}^{n+1}-\bu_{j,h}^n\|+\frac12\|\nabla\cdot\bu_{j,h}^{'n}\|_{L^\infty}\|\bu_{j,h}^{n+1}\|\|\bu_{j,h}^{n+1}-\bu_{j,h}^n\|\nonumber\\&\le \||\bu_{j,h}^{'n}|\nabla
\bu_{j,h}^{n+1}\|\|\bu_{j,h}^{n+1}-\bu_{j,h}^n\|+C\|\nabla\cdot\bu_{j,h}^{'n}\|_{L^\infty}\|\nabla\bu_{j,h}^{n+1}\|\|\bu_{j,h}^{n+1}-\bu_{j,h}^n\|\nonumber\\   &\le\frac{\alpha_j}{4}\|\nabla\bu_{j,h}^{n+1}\|^2+\|l_h^{n}\nabla 			\bu_{j,h}^{n+1}\|\|\bu_{j,h}^{n+1}-\bu_{j,h}^n\|+\frac{C}{\alpha_j}\|\nabla\cdot\bu_{j,h}^{'n}\|_{L^\infty}^2\|\bu_{j,h}^{n+1}-\bu_{j,h}^n\|^2\nonumber\\&\le\frac{\alpha_j}{4}\|\nabla\bu_{j,h}^{n+1}\|^2+\Delta t\|l_h^{n}\nabla 			\bu_{j,h}^{n+1}\|^2+\left(\frac{1}{4\Delta t}+\frac{C}{\alpha_j}\|\nabla\cdot\bu_{j,h}^{'n}\|_{L^\infty}^2\right)\|\bu_{j,h}^{n+1}-\bu_{j,h}^n\|^2.
  \end{align*}
  
Use of H\"older's and Young's inequalities, we have
\begin{align*}
    -(\nu_j^{'}\nabla \bu_{j,h}^{n},\nabla\bu_{j,h}^{n+1})\le\|\nu_j^{'}\|_{\infty}\|\nabla \bu_{j,h}^{n}\|\|\nabla\bu_{j,h}^{n+1}\|\le\frac{\|\nu_j^{'}\|_{\infty}}{2}\|\nabla \bu_{j,h}^{n}\|^2+\frac{\|\nu_j^{'}\|_{\infty}}{2}\|\nabla \bu_{j,h}^{n+1}\|^2.
\end{align*}

Using the above bounds, and reducing the equation \eqref{before-nl-bound-new}, becomes
	\begin{align}
		&\frac{1}{2\Delta t}\Big(\|\bu_{j,h}^{n+1}\|^2-\|\bu_{j,h}^n\|^2\Big)+\left(\frac{1}{4\Delta t}-\frac{C}{\alpha_j}\|\nabla\cdot\bu_{j,h}^{'n}\|_{L^\infty}^2\right)\|\bu_{j,h}^{n+1}-\bu_{j,h}^n\|^2+\frac{\Bar{\nu}_{\min}}{2}\|\nabla \bu_{j,h}^{n+1}\|^2\nonumber\\&+\gamma\|\nabla\cdot\bu_{j,h}^{n+1}\|^2+(2\mu-1)\Delta t\|l_h^n\nabla \bu_{j,h}^{n+1}\|^2\le \frac{1}{\alpha_j}\|\bif_{j}(t_{n+1})\|_{-1}^2+\frac{\|\nu_j^{'}\|_{\infty}}{2}\|\nabla \bu_{j,h}^{n}\|^2.\label{before-young}
	\end{align}
    
 Choose the calibration constant $\mu\ge\frac12$, time-step size given in \eqref{time-step-BE},
 drop non-negative term from left-hand-side, and rearrange
 \begin{align}
		\frac{1}{2\Delta t}\Big(\|\bu_{j,h}^{n+1}\|^2-\|\bu_{j,h}^n\|^2\Big)+\frac{\Bar{\nu}_{\min}}{2}\left(\|\nabla \bu_{j,h}^{n+1}\|^2-\|\nabla \bu_{j,h}^{n}\|^2\right)\nonumber\\+\frac{\alpha_j}{2}\|\nabla \bu_{j,h}^{n}\|^2+\gamma\|\nabla\cdot\bu_{j,h}^{n+1}\|^2\le \frac{1}{\alpha_j}\|\bif_{j}(t_{n+1})\|_{-1}^2.
	\end{align}
Multiply both sides by $2\Delta t$, and sum over the time-steps $n=0,1,\cdots,M-1$, which will finish the proof.
\end{proof}

\begin{theorem}\label{be-convergence} (Convergence of the BE-EEV scheme) Suppose $(\bu_j,p_j)$ satisfying \eqref{gov1-rans}-\eqref{Prandtl- mixing-length} and the following regularity assumptions for $m=\max\{3,k+1\}$
\begin{align*}
   \bu_j\in L^\infty(0,T; \bH^{m}(\cD)),
   \bu_{j,t}\in L^\infty(0,T;\bH^2(\cD)),
   \bu_{j,tt}\in L^\infty(0,T;\bL^2(\cD))
\end{align*}
with $k\ge 2$, then the ensemble solution of the Algorithm \ref{coupled-alg-com} converges to the true ensemble solution: For $\alpha_j>0$ and $\mu>\frac12$, if $$\Delta t\le\min_{\substack{1\le j\le J \\ 1\le n\le M}}\frac{C\alpha_j}{\|\nabla\cdot\bu_{j,h}^{'n}\|^2_{L^\infty}}$$ then, the following holds
\begin{align}
    \|<\bu>(T)-<\bu_h>^M\|^2+\alpha_{\min}\Delta t\sum_{n=1}^M\Big\|\nabla\Big(\hspace{-1.1mm}<\bu>(t_n)-<\bu_h>^n\hspace{-1.1mm}\Big)\Big\|^2\le C\Big(h^{2k}+\Delta t^2\Big).\label{convergence-error}
\end{align}
\end{theorem}
\begin{proof}
	We start our proof by forming the error equation. Testing \eqref{gov1-rans}-\eqref{gov2-rans} at the time level $t_{n+1}$, the continuous variational formulations can be written
	as
    
	\begin{align}
		\bigg(&\frac{\bu_{j}(t_{n+1})-\bu_{j}(t_n)}{\Delta t},\bv_h\bigg)+\Big(\Bar{\nu}\nabla \bu_{j}(t_{n+1}),\nabla\bv_h\Big)+\Big(\bu_j(t_{n+1})\cdot\nabla \bu_{j}(t_{n+1}),\bv_h\Big)\nonumber\\&+\left( 2\mu\Delta t(l^n)^2\nabla \bu_{j}(t_{n+1}),\nabla\bv_h\right)+\gamma\left(\nabla\cdot\bu_j(t_{n+1}),\nabla\cdot\bv_h\right)= \Big(\bif_{j}(t_{n+1}),\bv_h\Big)-\Big(\nu_j^{'}\nabla \bu_{j}(t_{n+1}),\nabla\bv_h\Big)\nonumber\\&+\left( 2\mu\Delta t \{(l^n)^2-(l^{n+1})^2\}\nabla \bu_{j}(t_{n+1}),\nabla\bv_h\right)-\bigg(\bu_{j,t}-\frac{\bu_{j}(t_{n+1})-\bu_{j}(t_n)}{\Delta t},\bv_h\bigg),\forall \bv_h\in \bV_h.\label{before-error-eqn-regular-new}
	\end{align}
    
	Denote $\be_j^n:=\bu_j(t_n)-\bu_{j,h}^n$. Set $\bchi_h=\bv_h\in \bV_h$ in \eqref{couple-eqn-1-new}, and then subtract \eqref{couple-eqn-1-new} from \eqref{before-error-eqn-regular-new}, to get
    
	\begin{align}
		&\bigg(\frac{\be_{j}^{n+1}-\be_{j}^n}{\Delta t},\bv_h\bigg)+\big(\Bar{\nu}\nabla \be_{j}^{n+1},\nabla\bv_h\big)+\left(\nu_j^{'}\nabla\be_{j}^n,\nabla\bv_h\right)+b^*\left(<\be>^n,\bu_j(t_{n+1})-\bu_{j}(t_n),\bv_h\right)\nonumber\\&+b^*\big(\hspace{-1.4mm}<\bu_h>^n,\be_j^{n+1},\bv_h\big)+b^*\left(\bu_{j,h}^{'n}, \be_{j}^n,\bv_h\right)+b^*\big(\be_j^n,\bu_{j}(t_n),\bv_h\big)+\left(2\mu\Delta t(l_h^n)^2\nabla \be_{j}^{n+1},\nabla \bv_h\right)\nonumber\\&+\left(2\mu\Delta t\{(l^n)^2-(l_h^n)^2\}\nabla \bu_{j}(t_{n+1}),\nabla \bv_h\right) +\gamma\left(\nabla\cdot\be_j^{n+1},\nabla\cdot\bv_h\right)=-G(t,\bu_j,\bv_h),
	\end{align}	
    
    where 
    
    \begin{align}
&G(t,\bu_j,\bv_h):=\bigg(\bu_{j,t}-\frac{\bu_{j}(t_{n+1})-\bu_{j}(t_n)}{\Delta t},\bv_h\bigg)+\left((\bu_j(t_{n+1})-\bu_j(t_{n}))\cdot\nabla \bu_{j}(t_{n+1}),\bv_h\right)\nonumber\\&+\big((\bu_j(t_n)-<\bu>(t_n))\cdot\nabla( \bu_{j}(t_{n+1})- \bu_{j}(t_{n})),\bv_h\big)+\left(\nu_j^{'}\nabla(\bu_{j}(t_{n+1})-\bu_{j}(t_n)),\nabla\bv_h\right)\nonumber\\&+\left(2\mu\Delta t \{(l^{n+1})^2-(l^n)^2\}\nabla \bu_{j}(t_{n+1}),\nabla\bv_h\right).
	\end{align}
    
Define $<\be>^n:=\frac{1}{J}\sum\limits_{j=1}^J\be_j^n$, and reduce as\begin{align*}
        &(l^n)^2-(l_h^n)^2=\sum_{j=1}^J\left(|\bu_{j}^{'}(t_n)|^2-|\bu_{j,h}^{'n}|^2\right)=\sum_{j=1}^J(\bu_{j}^{'}(t_n)-\bu_{j,h}^{'n})\cdot(\bu_{j}^{'}(t_n)+\bu_{j,h}^{'n})\\&=\sum_{j=1}^J\Big(\bu_j(t_n)-\frac{1}{J}\sum_{i=1}^J\bu_i(t_n)-\bu_{j,h}^n+\frac{1}{J}\sum_{i=1}^J\bu_{i,h}^n\Big)\cdot(\bu_{j}^{'}(t_n)+\bu_{j,h}^{'n})\\&=\sum_{j=1}^J(\be_j^n-<\be>^n)\cdot (\bu_{j}^{'}(t_n)+\bu_{j,h}^{'n}).
    \end{align*}
    
	Now, we decompose the error as the interpolation error and approximation term:
	\begin{align*}
		\be_{j}^n:& = \bu_j(t_n)-\bu_{j,h}^n=(\bu_j(t_n)-\tilde{\bu}_j^n)-(\bu_{j,h}^n-\tilde{\bu}_j^n):=\bfeta_{j}^n-\bphi_{j,h}^n,
	\end{align*}where $\tilde{\bu}_j^n: =P_{\bV_h}^{L^2}(\bu_j(t_n))\in \bV_h$ is the $L^2$ projections of $\bu_j(t_n)$ into $\bV_h$. Note that $\bphi_{j,h}^n\in\bV_h$ and $(\bfeta_{j}^n,\bv_{h})=0\hspace{2mm} \forall \bv_{h}\in \bV_h,$  
	we then have
    
	\begin{align}
&\bigg(\frac{\bphi_{j,h}^{n+1}-\bphi_{j,h}^n}{\Delta t},\bv_h\bigg)+\big(\Bar{\nu}\nabla \bphi_{j,h}^{n+1},\nabla\bv_h\big)+\left(\nu_j^{'}\nabla\bphi_{j,h}^n,\nabla\bv_h\right)+b^*\big(\hspace{-1.4mm}<\bu_h>^n,\bphi_{j,h}^{n+1},\bv_h\big)+b^*\left(\bu_{j,h}^{'n}, \bphi_{j,h}^n,\bv_h\right)\nonumber\\&+b^*\left(<\bphi_h>^n,\bu_j(t_{n+1})-\bu_{j}(t_n),\bv_h\right)+b^*\big(\bphi_{j,h}^n,\bu_{j}(t_n),\bv_h\big)+\gamma\left(\nabla\cdot\bphi_{j,h}^{n+1},\nabla\cdot\bv_h\right)\nonumber\\&+2\mu\Delta t\sum_{i=1}^J\left((\bphi_{i,h}^n-<\bphi_h>^n)\cdot (\bu_{i}^{'}(t_n)+\bu_{i,h}^{'n})\nabla \bu_{j}(t_{n+1}),\nabla \bv_h\right)+\left(2\mu\Delta t(l_h^n)^2\nabla \bphi_{j,h}^{n+1},\nabla \bv_h\right)\nonumber\\&=2\mu\Delta t\sum_{i=1}^J\left((\bfeta_i^n-<\bfeta>^n)\cdot (\bu_{i}^{'}(t_n)+\bu_{i,h}^{'n})\nabla \bu_{j}(t_{n+1}),\nabla \bv_h\right) +b^*\big(\hspace{-1.4mm}<\bu_h>^n,\bfeta_j^{n+1},\bv_h\big)\nonumber\\&+b^*\left(<\bfeta>^n,\bu_j(t_{n+1})-\bu_{j}(t_n),\bv_h\right)+b^*\left(\bu_{j,h}^{'n}, \bfeta_{j}^n,\bv_h\right)+b^*\big(\bfeta_j^n,\bu_{j}(t_n),\bv_h\big)+\big(\Bar{\nu}\nabla \bfeta_{j}^{n+1},\nabla\bv_h\big)\nonumber\\&+\left(\nu_j^{'}\nabla\bfeta_{j}^n,\nabla\bv_h\right)+\left(2\mu\Delta t(l_h^n)^2\nabla \bfeta_{j}^{n+1},\nabla \bv_h\right)+\gamma\left(\nabla\cdot\bfeta_{j}^{n+1},\nabla\cdot\bv_h\right)+G(t,\bu_j,\bv_h).\label{phi-equn-new}
	\end{align}
    
	Choose $\bv_h=\bphi_{j,h}^{n+1}$, use the polarization identity in \eqref{phi-equn-new}, and rearrange, to get
    
	\begin{align}
		&\frac{1}{2\Delta t}\lp\|\bphi_{j,h}^{n+1}\|^2-\|\bphi_{j,h}^{n}\|^2+\|\bphi_{j,h}^{n+1}-\bphi_{j,h}^{n}\|^2\rp+\big\|\Bar{\nu}^{\frac{1}{2}}\nabla\bphi_{j,h}^{n+1}\|^2+2\mu\Delta t\|l_h^n\nabla \bphi_{j,h}^{n+1}\|^2+\gamma\|\nabla\cdot\bphi_{j,h}^{n+1}\|^2\nonumber\\&+2\mu\Delta t\sum_{i=1}^J\left((\bphi_{i,h}^n-<\bphi_h>^n)\cdot (\bu_{i}^{'}(t_n)+\bu_{i,h}^{'n})\nabla \bu_{j}(t_{n+1}),\nabla \bphi_{j,h}^{n+1}\right)=-\left(\nu_j^{'}\nabla\bphi_{j,h}^n,\nabla\bphi_{j,h}^{n+1}\right)\nonumber\\&+\big(\Bar{\nu}\nabla \bfeta_{j}^{n+1},\nabla\bphi_{j,h}^{n+1}\big)+\left(\nu_j^{'}\nabla\bfeta_{j}^n,\nabla\bphi_{j,h}^{n+1}\right)-b^*\left(\bu_{j,h}^{'n}, \bphi_{j,h}^n,\bphi_{j,h}^{n+1}\right)-b^*\big(\bphi_{j,h}^n,\bu_{j}(t_n),\bphi_{j,h}^{n+1}\big)\nonumber\\&-b^*\left(<\bphi_h>^n,\bu_j(t_{n+1})-\bu_{j}(t_n),\bphi_{j,h}^{n+1}\right)+b^*\left(<\bfeta>^n,\bu_j(t_{n+1})-\bu_{j}(t_n),\bphi_{j,h}^{n+1}\right)\nonumber\\&+b^*\big(\hspace{-1.4mm}<\bu_h>^n,\bfeta_j^{n+1},\bphi_{j,h}^{n+1}\big)+b^*\left(\bu_{j,h}^{'n},\bfeta_{j}^n,\bphi_{j,h}^{n+1}\right)+b^*\big(\bfeta_j^n,\bu_{j}(t_n),\bphi_{j,h}^{n+1}\big)\nonumber\\&+\gamma\left(\nabla\cdot\bfeta_{j}^{n+1},\nabla\cdot\bphi_{j,h}^{n+1}\right)+\left(2\mu\Delta t(l_h^n)^2\nabla \bfeta_{j}^{n+1},\nabla\bphi_{j,h}^{n+1}\right)\nonumber\\&+2\mu\Delta t\sum_{i=1}^J\left((\bfeta_i^n-<\bfeta>^n)\cdot (\bu_{i}^{'}(t_n)+\bu_{i,h}^{'n})\nabla \bu_{j}(t_{n+1}),\nabla \bphi_{j,h}^{n+1}\right)+G\left(t,\bu_j,\bphi_{j,h}^{n+1}\right).\label{phibd-new}
	\end{align}
    
	Note that $b^*\big(\hspace{-1.4mm}<\bu_h>^n,\bphi_{j,h}^{n+1},\bphi_{j,h}^{n+1}\big)=0$. Now, turn our attention to finding bounds on the right side terms of \eqref{phibd-new}.
    Apply H\"older's and Young’s inequalities to obtain the following bounds\begin{align*}
		-\left(\nu_j^{'}\nabla\bphi_{j,h}^n,\nabla\bphi_{j,h}^{n+1}\right)\le\|\nu_j^{'}\|_\infty\|\nabla\bphi_{j,h}^n\|\|\nabla\bphi_{j,h}^{n+1}\|\le\frac{\|\nu_j^{'}\|_\infty}{2}\|\nabla\bphi_{j,h}^{n+1}\|^2+\frac{\|\nu_j^{'}\|_\infty}{2}\|\nabla\bphi_{j,h}^n\|^2,\\
		\big(\Bar{\nu}\nabla \bfeta_{j}^{n+1},\nabla\bphi_{j,h}^{n+1}\big)\le\|\Bar{\nu}\|_\infty\|\nabla\bfeta_{j}^{n+1}\|\|\nabla\bphi_{j,h}^{n+1}\|\le\frac{\alpha_j}{24}\|\nabla\bphi_{j,h}^{n+1}\|^2+\frac{6\|\Bar{\nu}\|_\infty^2}{\alpha_j}\|\nabla \bfeta_{j}^{n+1}\|^2,\\\left(\nu_j^{'}\nabla\bfeta_{j}^n,\nabla\bphi_{j,h}^{n+1}\right)\le\|\nu_j^{'}\|_\infty\|\bfeta_{j}^n\|\|\nabla\bphi_{j,h}^{n+1}\|\le\frac{\alpha_j}{24}\|\nabla\bphi_{j,h}^{n+1}\|^2+\frac{6\|\nu_j^{'}\|_\infty^2}{\alpha_j}\|\nabla \bfeta_{j}^{n+1}\|^2.
	\end{align*}
    
	For the first trilinear form, we rearrange, apply the identity \eqref{trilinear-identitiy}, Cauchy-Schwarz, H\"older's, Poincar\'e,  \eqref{basic-ineq}, and Young's inequalities, to have
	\begin{align*}
		&b^*\lp \bu_{j,h}^{'n},\bphi_{j,h} ^n,\bphi_{j,h}^{n+1}\rp=b^*\lp \bu_{j,h}^{'n},\bphi_{j,h} ^{n+1},\bphi_{j,h}^{n+1}-\bphi_{j,h}^{n}\rp\nonumber\\&=\lp \bu_{j,h}^{'n}\cdot\nabla\bphi_{j,h} ^{n+1},\bphi_{j,h}^{n+1}-\bphi_{j,h}^{n}\rp+\frac12\lp\nabla\cdot \bu_{j,h}^{'n},\bphi_{j,h} ^{n+1}\cdot\big(\bphi_{j,h}^{n+1}-\bphi_{j,h}^{n}\big)\rp\nonumber\\&\le\|\bu_{j,h}^{'n}\cdot \nabla\bphi_{j,h} ^{n+1}\|\|\bphi_{j,h}^{n+1}-\bphi_{j,h}^{n}\|+\frac12\|\nabla\cdot \bu_{j,h}^{'n}\|_{L^\infty}\|\bphi_{j,h} ^{n+1}\|\|\bphi_{j,h}^{n+1}-\bphi_{j,h}^{n}\|\nonumber\\&\le \|l_h^n \nabla\bphi_{j,h} ^{n+1}\|\|\bphi_{j,h}^{n+1}-\bphi_{j,h}^{n}\|+C\|\nabla\cdot \bu_{j,h}^{'n}\|_{L^\infty}\|\nabla\bphi_{j,h} ^{n+1}\|\|\bphi_{j,h}^{n+1}-\bphi_{j,h}^{n}\|\nonumber\\&\le\frac{\alpha_j}{24}\|\nabla\bphi_{j,h}^{n+1}\|^2+\Delta t\|l_h^n \nabla\bphi_{j,h} ^{n+1}\|^2+\left(\frac{1}{4\Delta t}+\frac{C}{\alpha_j}\|\nabla\cdot\bu_{j,h}^{'n}\|_{L^\infty}^2\right)\|\bphi_{j,h}^{n+1}-\bphi_{j,h}^{n}\|^2.
	\end{align*}
    
For the second trilinear term, we use the bound in \eqref{nonlinearbound3}, Agmon’s \cite{Robinson2016Three-Dimensional} inequality, Sobolev embedding theorem, regularity assumption of the true solution, and Young’s inequalities, to reveal\begin{align*}
		-b^*\lp\bphi_{j,h}^n, \bu_j(t_n),\bphi_{j,h}^{n+1}\rp&\le C\|\bphi_{j,h}^n\|\big(\|\nabla\bu_j(t_n)\|_{L^3}+\|\bu_j(t_n)\|_{L^\infty}\big)\|\nabla\bphi_{j,h}^{n+1}\|\le C\|\bphi_{j,h}^n\|\|\nabla\bphi_{j,h}^{n+1}\|\\&\le \frac{\alpha_j}{24}\|\nabla\bphi_{j,h}^{n+1}\|^2+\frac{C}{\alpha_j}\|\bphi_{j,h}^n\|^2.
	\end{align*}
    
 For the third trilinear term, use the bound in \eqref{nonlinearbound3}, triangle inequality, Agmon’s inequality, Sobolev embedding theorem, regularity assumption of the true solution $\bu_j\in L^\infty(0,T;\bH^3(\cD))$, and Young’s inequality, to obtain\begin{align*}
-&b^*\Big(\hspace{-1mm}<\hspace{-1mm}\bphi_{h}\hspace{-1mm}>^n,\bu_j(t_{n+1})-\bu_j(t_{n}),\bphi_{j,h}^{n+1}\Big)\\&\le  C\|<\hspace{-1mm}\bphi_{h}\hspace{-1mm}>^n\|\big(\|\nabla\left(\bu_j(t_{n+1})-\bu_j(t_{n})\right)\|_{L^3}+\|\bu_j(t_{n+1})-\bu_j(t_{n})\|_{L^\infty}\big)\|\nabla \bphi_{j,h}^{n+1}\|^2  \\&\le    \frac{\alpha_j}{24}\|\nabla \bphi_{j,h}^{n+1}\|^2+\frac{C}{\alpha_j}\|\hspace{-1mm}<\hspace{-1mm}\bphi_{h}\hspace{-1mm}>^n\hspace{-1mm}\|^2.
 \end{align*}
 
For the fourth trilinear term, apply the bound in \eqref{nonlinearbound}, triangle inequality, regularity assumption of the true solution, and Young’s inequality, to obtain \begin{align*}	
 b^*\Big(\hspace{-1mm}<\hspace{-1mm}\bfeta\hspace{-1mm}>^n,\bu_j(t_{n+1})-\bu_j(t_{n}),\bphi_{j,h}^{n+1}\Big)&\le C\|\nabla \hspace{-1mm}<\hspace{-1mm}\bfeta\hspace{-1mm}>^n\hspace{-1mm}\|\|\nabla\left(\bu_j(t_{n+1})-\bu_j(t_{n})\right)\|\|\nabla \bphi_{j,h}^{n+1}\|\\&\le \frac{\alpha_j}{24}\|\nabla \bphi_{j,h}^{n+1}\|^2+\frac{C}{\alpha_j}\|\nabla\hspace{-1mm} <\hspace{-1mm}\bfeta\hspace{-1mm}>^n\hspace{-1mm}\|^2.
 \end{align*}
 
For the fifth, and sixth trilinear terms, apply Young’s inequalities with \eqref{nonlinearbound}, to obtain \begin{align*}	
b^*\lp<\hspace{-1mm}\bu_h\hspace{-1mm}>^n, \bfeta_{j}^{n+1},\bphi_{j,h}^{n+1}\rp&\le C\|\nabla\hspace{-1mm}<\hspace{-1mm}\bu_h\hspace{-1mm}>^n\hspace{-1mm}\|\|\nabla \bfeta_{j}^{n+1}\|\|\nabla\bphi_{j,h}^{n+1}\|\\&\le \frac{\alpha_j}{24}\|\nabla \bphi_{j,h}^{n+1}\|^2+\frac{C}{\alpha_j}\|\nabla\hspace{-1mm}<\hspace{-1mm}\bu_h\hspace{-1mm}>^n\hspace{-1mm}\|^2\|\nabla \bfeta_{j}^{n+1}\|^2,\\
	b^*\lp \bu_{j,h}^{'n},\bfeta_{j} ^n,\bphi_{j,h}^{n+1}\rp&\le C\|\nabla \bu_{j,h}^{'n}\|\|\nabla\bfeta_{j} ^n\|\|\nabla\bphi_{j,h}^{n+1}\|\le \frac{\alpha_j}{24}\|\nabla \bphi_{j,h}^{n+1}\|^2+\frac{C}{\alpha_j}\|\nabla \bu_{j,h}^{'n}\|^2\|\nabla\bfeta_{j} ^n\|^2.
	\end{align*} 
    
 For the seventh trilinear term, apply the bound in \eqref{nonlinearbound}, regularity assumption of the true solution, and Young’s inequality, to obtain\begin{align*}	
 b^*\lp\bfeta^n_{j}, \bu_j(t_n),\bphi_{j,h}^{n+1}\rp&\le C\|\nabla \bfeta^n_{j}\|\|\nabla \bu_j(t_n)\|\|\nabla\bphi_{j,h}^{n+1}\|\le \frac{\alpha_j}{24}\|\nabla \bphi_{j,h}^{n+1}\|^2+\frac{C}{\alpha_j}\|\nabla \bfeta^n_{j}\|^2.
	\end{align*} 
    
Apply Cauchy-Schwarz, and Young’s inequalities,
	to obtain\begin{align*}
\gamma\left(\nabla\cdot\bfeta_{j}^{n+1},\nabla\cdot\bphi_{j,h}^{n+1}\right)&\le\frac{\gamma}{2}\|\nabla\cdot\bphi_{j,h}^{n+1}\|^2+C\gamma\|\nabla\bfeta_{j}^{n+1}\|^2.
	\end{align*}
    
   Assume $\mu>1/2$, apply Cauchy-Schwarz, Young’s, and H\"older's  inequalities, definition \eqref{eddy-viscosity-new}, triangle inequality and Lemma \ref{lemma1},
	to have\begin{align}
		2\mu\Delta t\lp (l_h^{n})^2\nabla \bfeta_{j}^{n+1},\nabla\bphi_{j,h}^{n+1}\rp&= 2\mu\Delta t\left(l_h^{n}\nabla \bfeta_{j}^{n+1},l_h^{n}\nabla\bphi_{j,h}^{n+1}\right)\le 2\mu\Delta t\|l_h^{n}\nabla \bfeta_{j}^{n+1}\|\|l_h^{n}\nabla\bphi_{j,h}^{n+1}\|\nonumber\\&\le \frac{2\mu-1}{2}\Delta t\|l_h^{n}\nabla\bphi_{j,h}^{n+1}\|^2+\frac{2\mu^2}{2\mu-1}\Delta t\|l_h^{n}\nabla \bfeta_{j}^{n+1}\|^2\nonumber\\&\le \frac{2\mu-1}{2}\Delta t\|l_h^{n}\nabla\bphi_{j,h}^{n+1}\|^2+\frac{2\mu^2}{2\mu-1}\Delta t\|(l_h^{n})^2\|_{L^\infty}\|\nabla \bfeta_{j}^{n+1}\|^2
        \nonumber\\&\le \frac{2\mu-1}{2}\Delta t\|l_h^{n}\nabla\bphi_{j,h}^{n+1}\|^2+\frac{2\mu^2}{2\mu-1}\Delta t\sum_{i=1}^J\|\bu_{i,h}^{'n}\|_{L^\infty}^2\|\nabla \bfeta_{j}^{n+1}\|^2\nonumber\\&\le \frac{2\mu-1}{2}\Delta t\|l_h^{n}\nabla\bphi_{j,h}^{n+1}\|^2+\frac{C\Delta t}{2\mu-1}\|\nabla \bfeta_{j}^{n+1}\|^2.\label{prandtl-mix-bound-new}
	\end{align}
    
    Using Taylor’s series expansion, Cauchy-Schwarz, H\"older's and Young's inequalities, the last term of \eqref{phibd-new} is bounded above as\begin{align*}
		\left|G(t,\bu_j, \bphi_{j,h}^{n+1})\right|\le\frac{\alpha_j}{24}\|\nabla\bphi_{j,h}^{n+1}\|^2+C\Delta t^2\Big(\|\bu_{j,tt}(s_1^*)\|^2+\|\nabla \bu_{j,t}(s_2^*)\|^2+\|\nabla \bu_{j,t}(s_3^{*})\|^2\|\nabla \bu_j(t_{n+1})\|^2\\+\|\bu_{j,t}(s_4^*)\|^2_{L^\infty}\|\nabla\bu_j(t_{n+1})\|^2+\|\nabla \big(\bu_j(t_n)-<\hspace{-1mm}\bu(t_n)\hspace{-1mm}>\hspace{-1mm}\big)\|^2\|\nabla \bu_{j,t}(s_5^{*})\|^2\Big),
	\end{align*} for some $s_i^*\in [t^n,t^{n+1}]$, $i=\overline{1,5}$. 
    
    Use H\"older's and triangle inequalities, regularity assumption of the true solution, Lemma \ref{lemma1}, Young's inequality and \begin{align}
\left(\sum\limits_{i=1}^J\|\ba_i\|\right)^2\le J\sum\limits_{i=1}^J\|\ba_i\|^2,\label{alg-ineq}
\end{align} to have
    \begin{align*}
        &2\mu\Delta t\sum_{i=1}^J\left((\bfeta_i^n-<\bfeta>^n)\cdot (\bu_{i}^{'}(t_n)+\bu_{i,h}^{'n})\nabla \bu_{j}(t_{n+1}),\nabla \bphi_{j,h}^{n+1}\right)\\&\le 2\mu\Delta t\sum_{i=1}^J\|\bfeta_i^n-<\bfeta>^n\|\|(\bu_{i}^{'}(t_n)+\bu_{i,h}^{'n})\cdot\nabla \bu_{j}(t_{n+1})\|_{L^\infty}\|\nabla \bphi_{j,h}^{n+1}\|\\&\le 2\mu\Delta t\sum_{i=1}^J\|\bfeta_i^n-<\bfeta>^n\|\|\bu_{i}^{'}(t_n)+\bu_{i,h}^{'n}\|_{L^\infty}\|\nabla \bu_{j}(t_{n+1})\|_{L^\infty}\|\nabla \bphi_{j,h}^{n+1}\|\\&\le CK^*\mu\Delta t\|\nabla\bphi_{j,h}^{n+1}\|\sum_{i=1}^J\|\bfeta_i^n\|\le\frac{\alpha_j}{24}\|\nabla\bphi_{j,h}^{n+1}\|^2+\frac{C\mu^2\Delta t^2}{\alpha_j}\sum_{i=1}^J\|\bfeta_i^n\|^2.
    \end{align*} 
    
    Similarly \begin{align*}
	    2\mu\Delta t\sum_{i=1}^J\left((\bphi_{i,h}^n-<\bphi_h>^n)\cdot (\bu_{i}^{'}(t_n)+\bu_{i,h}^{'n})\nabla \bu_{j}(t_{n+1}),\nabla \bphi_{j,h}^{n+1}\right)\\\le\frac{\alpha_j}{24}\|\nabla\bphi_{j,h}^{n+1}\|^2+\frac{C\mu^2\Delta t^2}{\alpha_j}\sum_{i=1}^J\|\bphi_{i,h}^n\|^2.
	\end{align*}
    
	Using these estimates in \eqref{phibd-new} and reducing, produces
	\begin{align}
		&\frac{1}{2\Delta t}\Big(\|\bphi_{j,h}^{n+1}\|^2-\|\bphi_{j,h}^{n}\|^2\Big)+\left(\frac{1}{4\Delta t}-\frac{C}{\alpha_j}\|\nabla\cdot\bu_{j,h}^{'n}\|_{L^\infty}^2\right)\|\bphi_{j,h}^{n+1}-\bphi_{j,h}^{n}\|^2+\frac{\Bar{\nu}_{\min}}{2}\|\nabla\bphi_{j,h}^{n+1}\|^2\nonumber\\&+\frac{\gamma}{2}\|\nabla\cdot\bphi_{j,h}^{n+1}\|^2+\frac{2\mu-1}{2}\Delta t\|l_h^n\nabla \bphi_{j,h}^{n+1}\|^2\le\frac{\|\nu_j^{'}\|_\infty}{2}\|\nabla\bphi_{j,h}^n\|^2+\frac{6}{\alpha_j}\left(\|\Bar{\nu}\|_\infty^2+\|\nu_j^{'}\|_\infty^2\right)\|\nabla \bfeta_{j}^{n+1}\|^2\nonumber\\&+\frac{C}{\alpha_j}\|\bphi_{j,h}^n\|^2+\frac{C}{\alpha_j}\|\hspace{-1mm}<\hspace{-1mm}\bphi_{h}\hspace{-1mm}>^n\hspace{-1mm}\|^2+\frac{C}{\alpha_j}\|\nabla\hspace{-1mm}<\hspace{-1mm}\bu_h\hspace{-1mm}>^n\hspace{-1mm}\|^2\|\nabla \bfeta_{j}^{n+1}\|^2+\frac{C}{\alpha_j}\|\nabla\hspace{-1mm} <\hspace{-1mm}\bfeta\hspace{-1mm}>^n\hspace{-1mm}\|^2\nonumber\\&+\frac{C}{\alpha_j}\|\nabla \bu_{j,h}^{'n}\|^2\|\nabla\bfeta_{j} ^n\|^2+\frac{C}{\alpha_j}\|\nabla \bfeta^n_{j}\|^2+\frac{C\mu^2\Delta t^2}{\alpha_j}\sum_{i=1}^J\left(\|\bfeta_i^n\|^2+\|\bphi_{i,h}^n\|^2\right)+\frac{C\Delta t}{2\mu-1}\|\nabla \bfeta_{j}^{n+1}\|^2\nonumber\\&+C\gamma\|\nabla\bfeta_{j}^{n+1}\|^2+C\Delta t^2\Big(\|\bu_{j,tt}(s_1^*)\|^2+\|\nabla \bu_{j,t}(s_2^*)\|^2+\|\nabla \bu_{j,t}(s_3^{*})\|^2\|\nabla \bu_j(t_{n+1})\|^2\nonumber\\&+\|\bu_{j,t}(s_4^*)\|^2_{\infty}\|\nabla\bu_j(t_{n+1})\|^2+\|\nabla \big(\bu_j(t_n)-<\hspace{-1mm}\bu(t_n)\hspace{-1mm}>\hspace{-1mm}\big)\|^2\|\nabla \bu_{j,t}(s_5^{*})\|^2\Big).
	\end{align}
    
	Assume $\mu>\frac12$, and choose time-step size $\Delta t\le\min\limits_{\substack{1\le j\le J \\ 1\le n\le M}}\frac{C\alpha_j}{\|\nabla\cdot\bu_{j,h}^{'n}\|^2_{L^\infty}}$, drop non-negative terms from left-hand-side, use \eqref{alg-ineq}, and rearrange
 \begin{align}
		\frac{1}{2\Delta t}\Big(&\|\bphi_{j,h}^{n+1}\|^2-\|\bphi_{j,h}^{n}\|^2\Big)+\frac{\Bar{\nu}_{\min}}{2}\left(\|\nabla\bphi_{j,h}^{n+1}\|^2-\|\nabla\bphi_{j,h}^{n}\|^2\right)+\frac{\alpha_j}{2}\|\nabla\bphi_{j,h}^{n}\|^2\nonumber\\&\le\frac{6}{\alpha_j}\left(\|\Bar{\nu}\|_\infty^2+\|\nu_j^{'}\|_\infty^2\right)\|\nabla \bfeta_{j}^{n+1}\|^2+\frac{C}{\alpha_j}\|\bphi_{j,h}^n\|^2+\frac{C}{\alpha_j}\sum_{i=1}^J\|\bphi_{i,h}^n\|^2+\frac{C}{\alpha_j}\|\nabla\hspace{-1mm}<\hspace{-1mm}\bu_h\hspace{-1mm}>^n\hspace{-1mm}\|^2\|\nabla \bfeta_{j}^{n+1}\|^2\nonumber\\&+\frac{C}{\alpha_j}\|\nabla\hspace{-1mm} <\hspace{-1mm}\bfeta\hspace{-1mm}>^n\hspace{-1mm}\|^2+\frac{C}{\alpha_j}\|\nabla \bu_{j,h}^{'n}\|^2\|\nabla\bfeta_{j} ^n\|^2+\frac{C}{\alpha_j}\|\nabla \bfeta^n_{j}\|^2+\frac{C\mu^2\Delta t^2}{\alpha_j}\sum_{i=1}^J\|\bfeta_i^n\|^2+\frac{C\Delta t}{2\mu-1}\|\nabla \bfeta_{j}^{n+1}\|^2\nonumber\\&+C\gamma\|\nabla\bfeta_{j}^{n+1}\|^2+C\Delta t^2\Big(\|\bu_{j,tt}(s_1^*)\|^2+\|\nabla \bu_{j,t}(s_2^*)\|^2+\|\nabla \bu_{j,t}(s_3^{*})\|^2\|\nabla \bu_j(t_{n+1})\|^2\nonumber\\&+\|\bu_{j,t}(s_4^*)\|^2_{\infty}\|\nabla\bu_j(t_{n+1})\|^2+\|\nabla \big(\bu_j(t_n)-<\hspace{-1mm}\bu(t_n)\hspace{-1mm}>\hspace{-1mm}\big)\|^2\|\nabla \bu_{j,t}(s_5^{*})\|^2\Big).
	\end{align}
    
 Multiplying both sides  by $2\Delta t$, sum over the time-steps $n=0,1,\cdots,M-1$, using $\|\bphi_{j,h}^0\|=\|\nabla\bphi_{j,h}^0\|=0$, $\Delta tM=T$, and using stability estimate, Agmon's inequality and regularity assumptions, to find
	\begin{align}   \|\bphi_{j,h}^{M}\|^2+\alpha_j\Delta t\sum_{n=1}^{M}\|\nabla\bphi_{j,h}^n\|^2\le C\Big(h^{2k}+\Delta t^2+\Delta t\sum_{n=1}^{M-1}\|\bphi_{j,h}^n\|^2+\Delta t\sum_{n=1}^{M-1}\sum_{i=1}^J\|\bphi_{i,h}^n\|^2\Big).\label{after-drop-non-negative-terms}
	\end{align}

Summing over $j=1,\cdots\hspace{-0.35mm},J$, and grouping, to get
	\begin{align}   \sum_{j=1}^J\|\bphi_{j,h}^{M}\|^2+\Delta t\sum_{n=1}^{M}\sum_{j=1}^J\alpha_j\|\nabla\bphi_{j,h}^n\|^2\le C\big(h^{2k}+\Delta t^2\big)+C\Delta t\sum_{n=1}^{M-1}\sum_{j=1}^J\|\bphi_{j,h}^n\|^2.
	\end{align}
    
Applying the discrete Gr\"onwall Lemma \ref{dgl}, we have\begin{align}   \sum_{j=1}^J\|\bphi_{j,h}^{M}\|^2+\Delta t\sum_{n=1}^{M}\sum_{j=1}^J\alpha_j\|\nabla\bphi_{j,h}^n\|^2\le C\Big(h^{2k}+\Delta t^2\Big).
	\end{align}
    
Now, using the triangle and Young's inequalities, we can write\begin{align}
\sum_{j=1}^J\|\be_{j}^{M}\|^2+\Delta t\sum_{n=1}^{M}\sum_{j=1}^J\alpha_j\|\nabla \be_{j}^{n}\|^2\le C\Big(h^{2k}+\Delta t^2\Big).\label{error-bounds-j-level}
	\end{align}
    
	Finally, again use the triangle and Young's inequalities to complete the proof.
\end{proof}

\subsection{Efficient linearized  BDF-2-EEV scheme} To obtain the second order BDF-2-EEV time-stepping scheme, we set the following in the Algorithm \ref{coupled-alg-com}
$$\beta=\frac{3}{2},\; \bsU_{j,h}^n=2\bu_{j,h}^n-\bu_{j,h}^{n-1},\; \bu_{j,h}^{'n}=2\bu_{j,h}^n-\bu_{j,h}^{n-1}-<\bu_h>^n,\;\;\text{and}\;\tilde{\bif}_{j}=\bif_{j}(t_{n+1})+\frac{2}{\Delta t}\bu_{j,h}^n-\frac{1}{2\Delta t}\bu_{j,h}^{n-1}.$$

For this $\bu_{j,h}^{'n}$, define $\tl_h^n:=\sum\limits_{j=1}^J|\bu_{j,h}^{'n}|^2$. We also define $
 \talpha_j := \Bar{\nu}_{\min} - 3 \|\nu^{'}_j\|_\infty > 0,\;\text{and}\;\talpha_{\min}:=\min\limits_{1\le j\le J}\talpha_j$. That is, $\frac{\|\nu^{'}_j\|_\infty}{\Bar{\nu}_{\min}}<\frac13.$
\begin{theorem} (Stability of the BDF-2-EEV scheme)\label{stability-bdf-2}
     Suppose $\bif_{\hspace{0.4mm}{j}} \in L^\infty \left(0,T,\bH^{-1}(\cD)\right)$, $\bu^0_{j,h}=\bu_j(\bx,0)\in \bH^2(\cD)\cap \bV$, $\bu^1_{j,h}=\bu_j(\bx,\Delta t) \in \bH^2(\cD)\cap \bV$ for $j=1,2,\cdots,J$, if \begin{align}
         \Delta t\le\min_{\substack{1\le j\le J \\ 1\le n\le M}}\frac{C\talpha_j}{\|\nabla\cdot\bu_{j,h}^{'n}\|^2_{L^\infty}},\label{time-step-size}
     \end{align} and choose $\mu \geq 1$ then, the solutions to the Algorithm \ref{coupled-alg-com} are stable,
\begin{align}
   &\|\bu_{j,h}^M\|^2+\|2\bu_{j,h}^M-\bu_{j,h}^{M-1}\|^2+2\talpha_{j} \Delta t \sum_{n=2}^{M}\|\nabla \bu_{j,h}^{n}\|^2 +4\gamma \Delta t \sum_{n=2}^{M}\|\nabla \cdot \bu_{j,h}^{n}\|^2 \leq \|\bu_{j,h}^1\|^2\nonumber \\
   &+\|2\bu_{j,h}^1-\bu_{j,h}^{0}\|^2+ 2\bar{\nu}_{\min}\Delta t \|\nabla \bu_{j,h}^{1}\|^2+ 2\Delta t\left(\bar{\nu}_{\min}-2\|\nu_{j}^{'}\|_{\infty}\right)\|\nabla \bu_{j,h}^{0}\|^2  +\frac{4\Delta t}{\talpha_j}\|\bif_{j}(t_{n+1})\|_{-1}^2.\label{stability-bdf2-statement}
\end{align}
\end{theorem}

\begin{proof}
    The proof follows by letting $\bchi_h = \bu^{n+1}_{j,h}$ in \eqref{couple-eqn-1-new}, and $q_h=p_{j,h}^{n+1}$ in \eqref{couple-incompressibility-new}, which gives
\begin{align}
    &\frac{1}{2\Delta t}\left(3\bu_{j,h}^{n+1}-4\bu_{j,h}^n+\bu_{j,h}^{n-1},\bu^{n+1}_{j,h}\right)+b^*\left(<\bu_h>^n, \bu_{j,h}^{n+1},\bu^{n+1}_{j,h}\right)+\left(\Bar{\nu}\nabla \bu_{j,h}^{n+1},\nabla\bu^{n+1}_{j,h}\right)\nonumber\\&+\gamma\left(\nabla\cdot\bu_{j,h}^{n+1},\nabla\cdot \bu^{n+1}_{j,h}\right)+\left(2\nu_T^h\nabla
    \bu_{j,h}^{n+1},\nabla \bu^{n+1}_{j,h}\right)= \left(\bif_{j}(t_{n+1}),\bu^{n+1}_{j,h}\right)\nonumber\\&-b^*(\bu_{j,h}^{'n}, 2\bu_{j,h}^{n}-\bu_{j,h}^{n-1},\bu^{n+1}_{j,h})-\left(\nu_j^{'}\nabla(2 \bu_{j,h}^{n}-\bu_{h}^{n-1}),\nabla\bu^{n+1}_{j,h}\right). \nonumber
\end{align}

Using $b^*\left(<\bu_h>^n, \bu_{j,h}^{n+1},\bu^{n+1}_{j,h}\right)=0$, \eqref{eddy-viscosity-new}, Cauchy-Schwarz and Young's inequalities, as well as the following algebraic identity
\begin{align}
\frac{1}{2}(3a - 4b + c)a = \frac{1}{4}\left[a^2 + \left(2a - b\right)^2\right] - \frac{1}{4}\left[b^2 + \left(2b - c\right)^2\right] + \frac{1}{4}\left(a - 2b + c\right)^2, \label{identity-bdf2}
\end{align}
we obtain,
\begin{align}
    &\frac{1}{4 \Delta t} \left( \|\bu^{n+1}_{j,h}\|^2 + \|2\bu^{n+1}_{j,h} - \bu^{n}_{j,h}\|^2 - \|\bu^{n}_{j,h}\|^2 - \|2\bu^{n}_{j,h} - \bu^{n-1}_{j,h}\|^2 + \|\bu^{n+1}_{j,h} - 2\bu^{n}_{j,h} + \bu^{n-1}_{j,h}\|^2\right) \notag \\ &+ \|\bar{\nu}^{\frac{1}{2}} \nabla \bu^{n+1}_{j,h}\|^2 + \gamma \| \nabla \cdot \bu^{n+1}_{j,h} \|^2 + 2\mu \Delta t \| \tl_h^n \nabla \bu^{n+1}_{j,h}\|^2 = \left(\bif_{j}(t_{n+1}), \bu^{n+1}_{j,h} \right) \notag \\ &- b^* \left( \bu^{'n}_{j,h} , 2\bu^n_{j,h} - \bu^{n-1}_{j,h} , \bu^{n+1}_{j,h} \right) - \left(\nu^{'}_j\nabla (2\bu^n_{j,h} - \bu^{n-1}_{j,h}) , \nabla \bu^{n+1}_{j,h}\right) . \label{Coupled-eq-1}
\end{align}

Using Cauchy-Schwarz, Young's, and H\"older's inequalities, yields
\begin{align}
    \Big(\bif_{j}(t_{n+1}), &\bu^{n+1}_{j,h} \Big) \leq \| \bif_{j}(t_{n+1}) \|_{-1} \| \nabla \bu^{n+1}_{j,h} \|
    \leq \frac{\talpha_j}{4} \| \nabla \bu^{n+1}_{j,h} \|^2 + \frac{1}{\talpha_j} \| \bif_{j}(t_{n+1}) \|^2_{-1}, \notag \\
    \Big(\nu^{'}_{j} \nabla (2\bu^n_{j,h} - \bu^{n-1}_{j,h}) , &\nabla \bu^{n+1}_{j,h}\Big)
    \leq 2 \|\nu^{'}_{j}\|_{\infty} \| \nabla \bu^n_{j,h}\|  \| \nabla \bu^{n+1}_{j,h} \| + \|\nu^{'}_{j}\|_{\infty} \| \nabla \bu^{n-1}_{j,h}\| \| \nabla \bu^{n+1}_{j,h} \| \notag \\
    &\leq \|\nu^{'}_j\|_{\infty} \| \nabla \bu^{n+1}_{j,h} \|^2 + \|\nu^{'}_j\|_{\infty} \|\nabla \bu^{n}_{j,h}\|^2 + \frac{1}{2} \|\nu^{'}_j\|_{\infty} \| \nabla \bu^{n+1}_{j,h} \|^2 + \frac{1}{2} \|\nu^{'}_j\|_{\infty} \|\nabla \bu^{n-1}_{j,h}\|^2 \notag
    \\
    &\leq \frac{3}{2}\|\nu^{'}_j \|_{\infty} \| \nabla \bu^{n+1}_{j,h} \|^2 + \|\nu^{'}_j\|_{\infty} \|\nabla \bu^{n}_{j,h}\|^2 + \frac{1}{2} \|\nu^{'}_j\|_{\infty} \|\nabla \bu^{n-1}_{j,h}\|^2 . \notag
\end{align}

Now, in the trilinear form using identity \eqref{trilinear-identitiy}, Cauchy-Schwarz, H\"older's, \eqref{basic-ineq}, Poincar\'e, and Young's inequalities, we get
\begin{align}
     &- b^* \left( \bu^{'n}_{j,h} , 2\bu^n_{j,h} - \bu^{n-1}_{j,h} , \bu^{n+1}_{j,h} \right)
     = b^* \left( \bu^{'n}_{j,h} , \bu^{n+1}_{j,h}, \bu_{j,h}^{n+1}-2\bu^n_{j,h} + \bu^{n-1}_{j,h} \right) \notag \\
     &= \left( \bu^{'n}_{j,h} \cdot \nabla \bu^{n+1}_{j,h}, (\bu_{j,h}^{n+1}-2\bu^n_{j,h} + \bu^{n-1}_{j,h}) \right) + \frac{1}{2} \left((\nabla \cdot \bu^{'n}_{j,h}) \bu^{n+1}_{j,h},\bu_{j,h}^{n+1}- 2\bu^n_{j,h} +\bu^{n-1}_{j,h} \right) \notag \\
     &\leq \| \bu^{'n}_{j,h} \cdot \nabla \bu^{n+1}_{j,h}\| \|\bu_{j,h}^{n+1}-2 \bu^n_{j,h}+ \bu^{n-1}_{j,h} \| + \frac{1}{2}\| \nabla \cdot \bu^{'n}_{j,h}\|_{L^\infty} \| \bu^{n+1}_{j,h} \| \|\bu_{j,h}^{n+1}-2 \bu^n_{j,h}+ \bu^{n-1}_{j,h} \| \notag \\
     &\leq \| |\bu^{'n}_{j,h}| \nabla \bu^{n+1}_{j,h} \| \|\bu_{j,h}^{n+1}-2 \bu^n_{j,h}+ \bu^{n-1}_{j,h} \| + C \| \nabla \cdot \bu^{'n}_{j,h}\|_{L^\infty} \| \nabla \bu^{n+1}_{j,h} \| \|\bu_{j,h}^{n+1}-2 \bu^n_{j,h}+ \bu^{n-1}_{j,h} \| \notag \\
     &\leq \|\tl_h^n \nabla\bu_{j,h}^{n+1}\|\|\bu_{j,h}^{n+1}-2 \bu^n_{j,h}+ \bu^{n-1}_{j,h} \|+ C \| \nabla \cdot \bu^{'n}_{j,h}\|_{L^\infty} \| \nabla \bu^{n+1}_{j,h} \| \|\bu_{j,h}^{n+1}-2 \bu^n_{j,h}+ \bu^{n-1}_{j,h} \| \notag \\
     &\leq \frac{\talpha_{j}}{4}\|\nabla\bu_{j,h}^{n+1}\|^2+2\Delta t \|\tl_h^n\nabla\bu_{j,h}^{n+1}\|^2+\left(\frac{1}{8\Delta t}+\frac{C}{\talpha_{j}}\|\nabla \cdot \bu_{j,h}^{'n}\|_{L^\infty}^{2}\right)\|\bu_{j,h}^{n+1}-2 \bu^n_{j,h}+ \bu^{n-1}_{j,h} \|^2. \nonumber
\end{align}

Using the above bounds, and rearranging, we have
\begin{align}
&\frac{1}{4 \Delta t} \left( \|\bu^{n+1}_{j,h}\|^2 + \|2\bu^{n+1}_{j,h} - \bu^{n}_{j,h}\|^2 - \|\bu^{n}_{j,h}\|^2 - \|2\bu^{n}_{j,h} - \bu^{n-1}_{j,h}\|^2 \right) \notag \\ & +\left(\frac{1}{8\Delta t}-\frac{C}{\talpha_{j}}\|\nabla \cdot \bu_{j,h}^{'n}\|_{L^\infty}^{2}\right)\|\bu_{j,h}^{n+1}-2 \bu^n_{j,h}+ \bu^{n-1}_{j,h} \|^2 + \frac{\bar{\nu}_{\min}}{2}  \| \nabla \bu^{n+1}_{j,h}\|^2 + \gamma \| \nabla \cdot \bu^{n+1}_{j,h} \|^2 \notag \\
& + 2\left(\mu-1\right)\Delta t \| \tl_h^n \nabla \bu^{n+1}_{j,h}\|^2 \leq \frac{1}{\talpha_j} \| \bif_{j}(t_{n+1}) \|^2_{-1} + \|\nu^{'}_j\|_{\infty} \| \nabla\bu^n_{j,h}\|^2 + \frac{1}{2} \|\nu^{'}_j\|_{\infty} \| \nabla\bu^{n-1}_{j,h}\|^2. \label{Coupled-eq-2}
\end{align}

Choose $\mu\ge 1$, time-step size given in \eqref{time-step-size}, drop the non-negative terms from the left-hand-side, and rearrange, yields
\begin{align}
    &\frac{1}{4\Delta t} \left( \|\bu^{n+1}_{j,h}\|^2- \|\bu^{n}_{j,h}\|^2 + \|2\bu^{n+1}_{j,h} - \bu^n_{j,h} \|^2 - \|2\bu^{n}_{j,h} - \bu^{n-1}_{j,h} \|^2 \right)  \notag\\&+ \frac{\bar{\nu}_{\min}}{2} \left(\| \nabla \bu^{n+1}_{j,h}\|^2 - \| \nabla \bu^{n}_{j,h}\|^2 \right) + \left(\frac{\bar{\nu}_{\min}}{2} - \|\nu^{'}_j\|_{\infty} \right) \left(\| \nabla \bu^{n}_{j,h}\|^2 - \| \nabla \bu^{n-1}_{j,h}\|^2 \right) \notag\\&+ \frac{\talpha_j}{2} \| \nabla \bu^{n-1}_{j,h} \|^2  + \gamma \| \nabla \cdot \bu^{n+1}_{j,h} \|^2 \leq \frac{1}{\talpha_j} \| \bif_{j}(t_{n+1}) \|^2_{-1}. 
\end{align}

Now, multiplying both sides by $4 \Delta t$, and summing over the time-steps $n=1,2,\cdots,M-1$, completes the proof.
\end{proof}

\begin{theorem} (Convergence of the BDF-2-EEV scheme)
Suppose $(\bu_j,p_j)$ satisfying \eqref{momentum}-\eqref{nse-initial} and the following regularity assumptions for $m=\max\{3,k+1\}$
\begin{align*}
   &\bu_j\in L^\infty(0,T; \bH^{m}(\cD)),\bu_{j,tt}\in L^\infty(0,T;\bH^2(\cD)), \bu_{j,ttt}\in L^\infty(0,T;\bL^2(\cD))
\end{align*}
with $k\ge 2$, then the ensemble solution of the BDF-2-EEV Algorithm \ref{coupled-alg-com} converges to the true ensemble solution: For $\talpha_j>0$ and $\mu> 1$, if $$\Delta t\le\min_{\substack{1\le j\le J \\ 1\le n\le M}}\frac{C\talpha_j}{\|\nabla\cdot\bu_{j,h}^{'n}\|^2_{L^\infty}}$$ then, the following holds
   \begin{align}
\|\bu_{j}(t_{M})-\bu_{j,h}^M\|^2+2\talpha_{\min}\Delta t\sum_{n=2}^{M}\sum_{j=1}^J\|\nabla\left(\bu_{j}(t_{n})-\bu_{j,h}^n\right)\|^2\le \frac{C}{\talpha_{\min}}(h^{2k}+\Delta t^4).
\end{align}
\end{theorem}
\begin{proof}
     Testing \eqref{gov1} with $\bv_{h}\in\bV_h$, at the time level $t_{n+1}$, the continuous variational formulations can be written
     
\begin{align}
\bigg(\frac{3\bu_j(t_{n+1})-4\bu_j(t_n)+\bu_j(t_{n-1})}{2\Delta t},\bv_{h}\bigg)+b^*\big(\bu_j(t_{n+1}), \bu_j(t_{n+1}),\bv_{h}\big)+\big(\bar{\nu}\nabla \bu_j(t_{n+1}), \nabla\bv_{h}\big)\nonumber\\+\big(\nu_j^{'}\nabla(2\bu_j(t_n)-\bu_j(t_{n-1})),\nabla\bv_{h}\big)+\left( 2\mu\Delta t(l^{n+1})^2\nabla \bu_{j}(t_{n+1}),\nabla\bv_h\right)+\gamma(\nabla\cdot\bu_{j}(t_{n+1}),\nabla\cdot \bv_{h})\nonumber\\=\big(\bif_{j}(t_{n+1}),\bv_{h}\big)-\Big(\nu^{'}_j\nabla\big(\bu_j(t_{n+1})-2\bu_j(t_n)+\bu_j(t_{n-1})\big),\nabla\bv_{h}\Big)\nonumber\\-\bigg(\bu_{j,t}(t_{n+1})-\frac{3\bu_j(t_{n+1})-4\bu_j(t_n)+\bu_j(t_{n-1})}{2\Delta t}, \bv_{h}\bigg). \label{conweakn1}
\end{align}

Set $\bchi_h=\bv_h\in \bV_h$ in \eqref{couple-eqn-1-new}, and then subtract \eqref{couple-eqn-1-new} from \eqref{conweakn1}, to get

\begin{align}
&\bigg(\frac{3\be^{n+1}_j-4\be^n_j+\be^{n-1}_j}{2\Delta t},\bv_{h}\bigg)+b^*\big(2\be_j^{n}-\be_j^{n-1}, \bu_j(t_{n+1}),\bv_{h}\big)+b^*\big(2\bu_{j,h}^{n}-\bu_{j,h}^{n-1}, \be_j^{n+1},\bv_{h}\big)\nonumber\\&+\gamma(\nabla\cdot\be_{j}^{n+1},\nabla\cdot \bv_{h})-b^*\big(\bu_{j,h}^{'n},\be_j^{n+1}-2\be_j^n+\be_j^{n-1},\bv_h\big)+\big(\bar{\nu}\nabla \be_j^{n+1}, \nabla\bv_{h}\big)+\big(\nu_j^{'}\nabla(2\be_j^n-\be_j^{n-1}),\nabla\bv_{h}\big)\nonumber\\&+ \left(2\mu\Delta t(\tl_h^n)^2\nabla \be_{j}^{n+1},\nabla \bv_h\right)+\left(2\mu\Delta t\{(l^{n+1})^2-(\tl_h^n)^2\}\nabla \bu_{j}(t_{n+1}),\nabla \bv_h\right)=-G(t,\bu_j,\bv_h), \label{error}
\end{align}
where
\begin{align}
G(t,\bu_j,\bv_h):=b^*\big(\bu_{j,h}^{'n},\bu_j(t_{n+1})-2\bu_j(t_n)+\bu_j(t_{n-1}),\bv_h\big)\nonumber\\+b^*\big(\bu_j(t_{n+1})-2\bu_j(t_n)+\bu_j(t_{n-1}),\bu_j(t_{n+1}),\bv_h\big)\nonumber\\+\Big(\nu^{'}_j\nabla\big(\bu_j(t_{n+1})-2\bu_j(t_n)+\bu_j(t_{n-1})\big),\nabla\bv_{h}\Big)\nonumber\\+\bigg(\bu_{j,t}(t_{n+1})-\frac{3\bu_j(t_{n+1})-4\bu_j(t_n)+\bu_j(t_{n-1})}{2\Delta t}, \bv_{h}\bigg).\label{G}
\end{align}

Define $<\be>^n:=\frac{1}{J}\sum\limits_{j=1}^J(2\be_j^n-\be_j^{n-1})$, and consider the following using Taylor series expansion 

\begin{align*}
        &\bu_{j}^{'}(t_{n+1})-\bu_{j,h}^{'n}=\bu_j(t_{n+1})-\frac{1}{J}\sum_{i=1}^J\bu_j(t_{n+1})-\big(2\bu_{j,h}^n-\bu_{j,h}^{n-1}\big)+\frac{1}{J}\sum_{i=1}^J\big(2\bu_{i,h}^n-\bu_{i,h}^{n-1}\big)\nonumber\\&=2\be_j^n-\be_j^{n-1}-<\be>^n+\bu_j(t_{n+1})-2\bu_j(t_n)+\bu_j(t_{n-1})-\frac{1}{J}\sum_{j=1}^J(\bu_j(t_{n+1})-2\bu_j(t_n)+\bu_j(t_{n-1}))\nonumber\\&=2\be_j^n-\be_j^{n-1}-<\be>^n+\Delta t^2\bu_{j,tt}(t_1^*)-\frac{1}{J}\Delta t^2\sum_{j=1}^J\bu_{j,tt}(t_1^*),\hspace{3mm}\text{for some}\hspace{2mm}t_1^*\in [t^{n-1},t^{n+1}].
    \end{align*} 
    
    Now, reduce $(l^{n+1})^2-(\tl_h^n)^2$ as below
    
    \begin{align}
        (l^{n+1})^2&-(\tl_h^n)^2=\sum_{j=1}^J\left(|\bu_{j}^{'}(t_{n+1})|^2-|\bu_{j,h}^{'n}|^2\right)=\sum_{j=1}^J(\bu_{j}^{'}(t_{n+1})-\bu_{j,h}^{'n})\cdot(\bu_{j}^{'}(t_{n+1})+\bu_{j,h}^{'n})\nonumber\\=&\sum_{j=1}^J\bigg(2\be_j^n-\be_j^{n-1}-<\be>^n+\Delta t^2\bu_{j,tt}(t_1^*)-\frac{1}{J}\Delta t^2\sum_{j=1}^J\bu_{j,tt}(t_1^*)\bigg)\cdot (\bu_{j}^{'}(t_{n+1})+\bu_{j,h}^{'n}).
    \end{align}
    
	Now, we decompose the error as the interpolation error and approximation term:
    
	\begin{align*}
		\be_{j}^n:& = \bu_j(t_n)-\bu_{j,h}^n=(\bu_j(t_n)-\tilde{\bu}_j^n)-(\bu_{j,h}^n-\tilde{\bu}_j^n):=\bfeta_{j}^n-\bphi_{j,h}^n,
	\end{align*}
	where $\tilde{\bu}_j^n: =P_{\bV_h}^{L^2}(\bu_j(t_n))\in \bV_h$ is the $L^2$ projections of $\bu_j(t_n)$ into $\bV_h$, 
	we then have

\begin{align}
\frac{1}{2\Delta t}\Big(3\bphi_{j,h}^{n+1}-4\bphi_{j,h}^n+\bphi_{j,h}^{n-1},\bv_{h}\Big)+\big(\Bar{\nu} \nabla \bphi_{j,h}^{n+1},\nabla\bv_{h}\big)
    +b^*\big(2\bu_{j,h}^n-\bu_{j,h}^{n-1},\bphi_{j,h}^{n+1},\bv_{h}\big)\nonumber\\+b^*\big(2\bphi_{j,h}^n-\bphi_{j,h}^{n-1},\bu_{j}(t_{n+1}),\bv_{h}\big)
   -b^*\big(\bu_{j,h}^{'n}, \bphi_{j,h}^{n+1}-2\bphi_{j,h}^{n}+\bphi_{j,h}^{n-1},\bv_{h}\big)
\nonumber\\+\gamma\big(\nabla\cdot\bphi_{j,h}^{n+1},\nabla\cdot\bv_h\big)+\big( \nu_j^{'}\nabla(2 \bphi_{j,h}^{n}-\bphi_{j,h}^{n-1}),\nabla\bv_{h}\big)+\left(2\mu\Delta t(\tl_h^n)^2\nabla \bphi_{j,h}^{n+1},\nabla \bv_h\right)\nonumber\\+2\mu\Delta t\sum_{i=1}^J\left((2\bphi_{i,h}^n-\bphi_{i,h}^{n-1}-<\bphi_h>^n)\cdot (\bu_{i}^{'}(t_{n+1})+\bu_{i,h}^{'n})\nabla \bu_{j}(t_{n+1}),\nabla \bv_h\right)\nonumber\\+2\mu\Delta t^3\sum_{i=1}^J\bigg(\big(\bu_{j,tt}(t_1^*)-\frac{1}{J}\sum_{j=1}^J\bu_{j,tt}(t_1^*)\big)\cdot (\bu_{i}^{'}(t_{n+1})+\bu_{i,h}^{'n})\nabla \bu_{j}(t_{n+1}),\nabla \bv_h\bigg)\nonumber\\=2\mu\Delta t\sum_{i=1}^J\left((2\bfeta_i^n-\bfeta_i^{n-1}-<\bfeta>^n)\cdot (\bu_{i}^{'}(t_{n+1})+\bu_{i,h}^{'n})\nabla \bu_{j}(t_{n+1}),\nabla \bv_h\right)
  -\big(\Bar{\nu} \nabla \bfeta_j^{n+1},\nabla\bv_{h}\big)\nonumber\\
  -\gamma\big(\nabla\cdot\bfeta_j^{n+1},\nabla\cdot\bv_h\big)-\big( \nu_j^{'}\nabla(2 \bfeta_{j}^{n}-\bfeta_{j}^{n-1}),\nabla\bv_{h}\big)-b^*\big(2\bu_{j,h}^n-\bu_{j,h}^{n-1},\bfeta_{j}^{n+1},\bv_{h}\big)\nonumber\\-b^*\big(2\bfeta_{j}^n-\bfeta_{j}^{n-1},\bu_{j}(t_{n+1}),\bv_{h}\big)
   +b^*\big(\bu_{j,h}^{'n}, \bfeta_{j}^{n+1}-2\bfeta_{j}^{n}+\bfeta_{j}^{n-1},\bv_{h}\big)\nonumber\\+\left(2\mu\Delta t(\tl_h^n)^2\nabla \bfeta_{j}^{n+1},\nabla \bv_h\right)
 -G(t,\bu_j,\bv_h).\label{before-phi-bdf2}
\end{align}

Now, choose $\bv_h=\bphi_{j,h}^{n+1}$, use the algebraic identity \eqref{identity-bdf2}, and rearrange to obtain

\begin{align}
    &\frac{1}{4\Delta t}\Big(\|\bphi_{j,h}^{n+1}\|^2
    +\|2\bphi_{j,h}^{n+1}-\bphi_{j,h}^{n}\|^2-\|\bphi_{j,h}^{n}\|^2-\|2\bphi_{j,h}^{n}-\bphi_{j,h}^{n-1}\|^2
+\|\bphi_{j,h}^{n+1}-2\bphi_{j,h}^{n}
    +\bphi_{j,h}^{n-1}\|^2\Big)\nonumber\\
    &+\|\Bar{\nu}^{\frac{1}{2}}\nabla\bphi_{j,h}^{n+1}\|^2
+\gamma\|\nabla\cdot\bphi_{j,h}^{n+1}\|^2+\big( \nu_j^{'}\nabla(2 \bphi_{j,h}^{n}-\bphi_{j,h}^{n-1}),\nabla\bphi_{j,h}^{n+1}\big)+2\mu\Delta t\|\tl_h^n\nabla \bphi_{j,h}^{n+1}\|^2\nonumber\\&+2\mu\Delta t\sum_{i=1}^J\left((2\bphi_{i,h}^n-\bphi_{i,h}^{n-1}-<\bphi_h>^n)\cdot (\bu_{i}^{'}(t_n)+\bu_{i,h}^{'n})\nabla \bu_{j}(t_{n+1}),\nabla \bphi_{j,h}^{n+1}\right)\nonumber\\&+2\mu\Delta t^3\sum_{i=1}^J\bigg(\big(\bu_{j,tt}(t_1^*)-\frac{1}{J}\sum_{j=1}^J\bu_{j,tt}(t_1^*)\big)\cdot (\bu_{i}^{'}(t_{n+1})+\bu_{i,h}^{'n})\nabla \bu_{j}(t_{n+1}),\nabla \bphi_{j,h}^{n+1}\bigg)\nonumber\\&=2\mu\Delta t\sum_{i=1}^J\left((2\bfeta_i^n-\bfeta_i^{n-1}-<\bfeta>^n)\cdot (\bu_{i}^{'}(t_n)+\bu_{i,h}^{'n})\nabla \bu_{j}(t_{n+1}),\nabla \bphi_{j,h}^{n+1}\right)
    \nonumber\\&-\big(\Bar{\nu} \nabla \bfeta_j^{n+1},\nabla\bphi_{j,h}^{n+1}\big)-\gamma\big(\nabla\cdot\bfeta_j^{n+1},\nabla\cdot\bphi_{j,h}^{n+1}\big)-b^*\big(2\bphi_{j,h}^n-\bphi_{j,h}^{n-1},\bu_{j}(t_{n+1}),\bphi_{j,h}^{n+1}\big)
    \nonumber\\&+b^*\big(\bu_{j,h}^{'n}, \bphi_{j,h}^{n+1}-2\bphi_{j,h}^{n}+\bphi_{j,h}^{n-1},\bphi_{j,h}^{n+1}\big)-b^*\big(2\bu_{j,h}^n-\bu_{j,h}^{n-1},\bfeta_{j}^{n+1},\bphi_{j,h}^{n+1}\big)\nonumber\\&
  -b^*\big(2\bfeta_{j}^n-\bfeta_{j}^{n-1},\bu_{j}(t_{n+1}),\bphi_{j,h}^{n+1}\big)
  +b^*\big(\bu_{j,h}^{'n}, \bfeta_{j}^{n+1}-2\bfeta_{j}^{n}+\bfeta_{j}^{n-1},\bphi_{j,h}^{n+1}\big)\nonumber\\
  &
  -\big( \nu_j^{'}\nabla(2 \bfeta_{j}^{n}-\bfeta_{j}^{n-1}),\nabla\bphi_{j,h}^{n+1}\big)+\left(2\mu\Delta t(\tl_h^n)^2\nabla \bfeta_{j}^{n+1},\nabla\bphi_{j,h}^{n+1}\right)-G(t,\bu_j,\bphi_{j,h}^{n+1}).\label{orthdecomp}
\end{align}

Using H\"older's, triangle and Young's inequalities,

\begin{align*}
    -\big( \nu_j^{'}\nabla(2 \bphi_{j,h}^{n}-\bphi_{j,h}^{n-1}),\nabla 
		\bphi_{j,h}^{n+1}\big)
  &\leq 2\|\nu_j^{'}\|_{\infty}\|\nabla 
		\bphi_{j,h}^{n}\|\|\nabla 
		\bphi_{j,h}^{n+1}\|+\|\nu_j^{'}\|_{\infty}\|\nabla 
		\bphi_{j,h}^{n-1}\|\|\nabla 
		\bphi_{j,h}^{n+1}\|\\
  &\leq \frac{3}{2}\|\nu_j^{'}\|_{\infty}\|\nabla 
		\bphi_{j,h}^{n+1}\|^2+\|\nu_j^{'}\|_{\infty}\|\nabla 
		\bphi_{j,h}^{n}\|^2+\frac{1}{2}\|\nu_j^{'}\|_{\infty}\|\nabla 
		\bphi_{j,h}^{n-1}\|^2,\\
          -\big( \nu_j^{'}\nabla(2 \bfeta_{j}^{n}-\bfeta_{j}^{n-1}),\nabla\bphi_{j,h}^{n+1}\big)&\le 2\|\nu_j^{'}\|_{\infty}\|\nabla 
		\bfeta_{j}^{n}\|\|\nabla 
		\bphi_{j}^{n+1}\|+\|\nu_j^{'}\|_{\infty}\|\nabla 
		\bfeta_{j}^{n-1}\|\|\nabla 
		\bphi_{j}^{n+1}\|\\
  &\leq \frac{\talpha_j}{22}\|\nabla 
		\bphi_{j}^{n+1}\|^2+\frac{44\|\nu_j^{'}\|_{\infty}^2}{\talpha_j}\|\nabla 
		\bfeta_{j}^{n}\|^2+\frac{11\|\nu_j^{'}\|_{\infty}}{\talpha_j}\|\nabla 
		\bfeta_{j}^{n-1}\|^2.
\end{align*}

Apply H\"older's and triangle inequalities, regularity assumption of the true solution, Lemma \ref{lemma1}, Young's and \eqref{alg-ineq} inequalities  to obtain 

\begin{align*}
	    -&2\mu\Delta t\sum_{i=1}^J\left((2\bphi_{i,h}^n-\bphi_{i,h}^{n-1}-<\bphi_h>^n)\cdot (\bu_{i}^{'}(t_n)+\bu_{i,h}^{'n})\nabla \bu_{j}(t_{n+1}),\nabla \bphi_{j,h}^{n+1}\right)\\\le &2\mu\Delta t\sum_{i=1}^J\|2\bphi_{i,h}^n-\bphi_{i,h}^{n-1}-<\bphi_h>^n\|\|(\bu_{i}^{'}(t_n)+\bu_{i,h}^{'n})\cdot\nabla \bu_{j}(t_{n+1})\|_{L^\infty}\|\nabla\bphi_{j,h}^{n+1}\|\\\le &2\mu\Delta t\sum_{i=1}^J\|2\bphi_{i,h}^n-\bphi_{i,h}^{n-1}-<\bphi_h>^n\|\|\bu_{i}^{'}(t_n)+\bu_{i,h}^{'n}\|_{L^\infty}\|\nabla \bu_{j}(t_{n+1})\|_{L^\infty}\|\nabla\bphi_{j,h}^{n+1}\|\\\le &CK^*\mu\Delta t\|\nabla\bphi_{j,h}^{n+1}\|\sum_{i=1}^J\left(\|\bphi_{i,h}^n\|+\|\bphi_{i,h}^{n-1}\|\right)\le\frac{\talpha_j}{22}\|\nabla\bphi_{j,h}^{n+1}\|^2+\frac{C\mu^2\Delta t^2}{\talpha_j}\sum_{i=1}^J\left(\|\bphi_{i,h}^n\|^2+\|\bphi_{i,h}^{n-1}\|^2\right).
	\end{align*}
    
    Similarly
    \begin{align*}
        2\mu\Delta t\sum_{i=1}^J\left((2\bfeta_i^n-\bfeta_i^{n-1}-<\bfeta>^n)\cdot (\bu_{i}^{'}(t_n)+\bu_{i,h}^{'n})\nabla \bu_{j}(t_{n+1}),\nabla \bphi_{j,h}^{n+1}\right)\\\le\frac{\talpha_j}{22}\|\nabla\bphi_{j,h}^{n+1}\|^2+\frac{C\mu^2\Delta t^2}{\talpha_j}\sum_{i=1}^J\left(\|\bfeta_i^n\|^2+\|\bfeta_{i}^{n-1}\|^2\right).
    \end{align*}
    
Applying H\"older's and Young's inequalities,

\begin{align*}
        -\big(\Bar{\nu}\nabla \bfeta_{j}^{n+1},\nabla\bphi_{j,h}^{n+1}\big)\le\|\Bar{\nu}\|_\infty\|\nabla\bfeta_{j}^{n+1}\|\|\nabla\bphi_{j,h}^{n+1}\|\le\frac{\talpha_j}{22}\|\nabla\bphi_{j,h}^{n+1}\|^2+\frac{11\|\Bar{\nu}\|_\infty^2}{2\talpha_j}\|\nabla \bfeta_{j}^{n+1}\|^2.
\end{align*}

Apply Cauchy-Schwarz, and Young’s inequalities,
	to obtain
    
	\begin{align*}
-\gamma\left(\nabla\cdot\bfeta_{j}^{n+1},\nabla\cdot\bphi_{j,h}^{n+1}\right)&\le\frac{\gamma}{2}\|\nabla\cdot\bphi_{j,h}^{n+1}\|^2+C\gamma\|\nabla\bfeta_{j}^{n+1}\|^2.
	\end{align*}
    
Using \eqref{nonlinearbound3}, regularity assumption on $\bu_{j}(t_{n+1})$, Young's, and triangle inequalities,

\begin{align}
    -b^*\big(2\bphi_{j,h}^n-\bphi_{j,h}^{n-1},\bu_{j}(t_{n+1}),\bphi_{j,h}^{n+1}\big)
    &\leq C\|2\bphi_{j,h}^n-\bphi_{j,h}^{n-1}\| \left(\|\nabla \bu_{j}(t_{n+1})\|_{L^3}+\|\bu_{j}(t_{n+1})\|_{L^\infty}\right) \|\nabla\bphi_{j,h}^{n+1}\|\nonumber\\
    &\leq C \|2\bphi_{j,h}^n-\bphi_{j,h}^{n-1}\|\|\nabla\bphi_{j,h}^{n+1}\|\nonumber\\
    &\leq \frac{\talpha_j}{22}\|\nabla\bphi_{j,h}^{n+1}\|^2+\frac{C}{\talpha_j}\|2\bphi_{j,h}^n-\bphi_{j,h}^{n-1}\|^2\nonumber\\
    &\leq \frac{\talpha_j}{22}\|\nabla\bphi_{j,h}^{n+1}\|^2+\frac{C}{\talpha_j}\Big(\|\bphi_{j,h}^n\|^2+\|\bphi_{j,h}^{n-1}\|^2\Big).\nonumber
\end{align}

Apply identity in \eqref{trilinear-identitiy}, Cauchy-Schwarz, H\"older's, \eqref{basic-ineq}, Poincar\'e, and Young's inequalities, we get

 \begin{align*}
        &b^*\big(\bu_{j,h}^{'n}, \bphi_{j,h}^{n+1}-2\bphi_{j,h}^{n}+\bphi_{j,h}^{n-1},\bphi_{j,h}^{n+1}\big)=-b^*\big(\bu_{j,h}^{'n},\bphi_{j,h}^{n+1}, \bphi_{j,h}^{n+1}-2\bphi_{j,h}^{n}+\bphi_{j,h}^{n-1}\big)\\&=-\left(\bu_{j,h}^{'n}\cdot\nabla\bphi_{j,h}^{n+1},\bphi_{j,h}^{n+1}-2\bphi_{j,h}^{n}+\bphi_{j,h}^{n-1}\right)-\frac12\left(\nabla\cdot\bu_{j,h}^{'n},\bphi_{j,h}^{n+1}\cdot\big(\bphi_{j,h}^{n+1}-2\bphi_{j,h}^{n}+\bphi_{j,h}^{n-1}\big)\right)\\&\le\|\bu_{j,h}^{'n}\cdot\nabla\bphi_{j,h}^{n+1}\|\|\bphi_{j,h}^{n+1}-2\bphi_{j,h}^{n}+\bphi_{j,h}^{n-1}\|+\frac12\|\nabla\cdot\bu_{j,h}^{'n}\|_{L^\infty}\|\bphi_{j,h}^{n+1}\|\|\bphi_{j,h}^{n+1}-2\bphi_{j,h}^{n}+\bphi_{j,h}^{n-1}\|\\&\le\|l_h^n\nabla\bphi_{j,h}^{n+1}\|\|\bphi_{j,h}^{n+1}-2\bphi_{j,h}^{n}+\bphi_{j,h}^{n-1}\|+C\|\nabla\cdot\bu_{j,h}^{'n}\|_{L^\infty}\|\nabla\bphi_{j,h}^{n+1}\|\|\bphi_{j,h}^{n+1}-2\bphi_{j,h}^{n}+\bphi_{j,h}^{n-1}\|\\&\le\frac{\talpha_j}{22}\|\bphi_{j,h}^{n+1}\|^2+ 2\Delta t\|\tl_h^n\nabla\bphi_{j,h}^{n+1}\|^2+\left(\frac{1}{8\Delta t}+\frac{C}{\talpha_j}\|\nabla\cdot\bu_{j,h}^{'n}\|_{L^\infty}^2\right)\|\bphi_{j,h}^{n+1}-2\bphi_{j,h}^{n}+\bphi_{j,h}^{n-1}\|^2.
    \end{align*}
    
Use the non-linear bound in \eqref{nonlinearbound} and Young's inequality to get

    \begin{align*}
        -b^*\big(2\bu_{j,h}^n-\bu_{j,h}^{n-1},\bfeta_{j}^{n+1},\bphi_{j,h}^{n+1}\big)&\le C\|\nabla(2\bu_{j,h}^n-\bu_{j,h}^{n-1})\|\|\nabla\bfeta_{j}^{n+1}\|\|\nabla\bphi_{j,h}^{n+1}\|\\&\le \frac{\talpha_j}{22}\|\nabla\bphi_{j,h}^{n+1}\|^2+\frac{C}{\talpha_j}\|\nabla(2\bu_{j,h}^n-\bu_{j,h}^{n-1})\|^2\|\nabla\bfeta_{j}^{n+1}\|^2,\\
        b^*\big(\bu_{j,h}^{'n}, \bfeta_{j}^{n+1}-2\bfeta_{j}^{n}+\bfeta_{j}^{n-1},\bphi_{j,h}^{n+1}\big)&\le C\|\nabla \bu_{j,h}^{'n}\|\|\nabla(\bfeta_{j}^{n+1}-2\bfeta_{j}^{n}+\bfeta_{j}^{n-1})\|\|\nabla \bphi_{j,h}^{n+1}\|\\&\le \frac{\talpha_j}{22}\|\nabla\bphi_{j,h}^{n+1}\|^2+\frac{C}{\talpha_j}\|\nabla \bu_{j,h}^{'n}\|^2\|\nabla(\bfeta_{j}^{n+1}-2\bfeta_{j}^{n}+\bfeta_{j}^{n-1})\|^2.
    \end{align*}
    
    Apply the non-linear bound in \eqref{nonlinearbound}, the regularity assumption, and Young's inequality to get
    \begin{align*}
        -b^*\big(2\bfeta_{j}^n-\bfeta_{j}^{n-1},\bu_{j}(t_{n+1}),\bphi_{j,h}^{n+1}\big)&\le C\|\nabla(2\bfeta_{j}^n-\bfeta_{j}^{n-1})\|\|\nabla \bu_{j}(t_{n+1})\|\|\nabla \bphi_{j,h}^{n+1}\|\\&\le C\|\nabla(2\bfeta_{j}^n-\bfeta_{j}^{n-1})\|\|\nabla \bphi_{j,h}^{n+1}\|\\&\le \frac{\talpha_j}{22}\|\nabla\bphi_{j,h}^{n+1}\|^2+\frac{C}{\talpha_j}\|\nabla(2\bfeta_{j}^n-\bfeta_{j}^{n-1})\|^2.
    \end{align*}

Use of H\"older's, and triangle inequalities, regularity assumption, Lemma \ref{lemma1} and Young's inequality, provides
\begin{align*}
    &2\mu\Delta t^3\sum_{i=1}^J\bigg(\big(\bu_{j,tt}(t_1^*)-\frac{1}{J}\sum_{j=1}^J\bu_{j,tt}(t_1^*)\big)\cdot (\bu_{i}^{'}(t_{n+1})+\bu_{i,h}^{'n})\nabla \bu_{j}(t_{n+1}),\nabla \bphi_{j,h}^{n+1}\bigg)\\&\le 2\mu\Delta t^3\sum_{i=1}^J\Big\|\bu_{j,tt}(t_1^*)-\frac{1}{J}\sum_{j=1}^J\bu_{j,tt}(t_1^*)\Big\|_{L^\infty}\|\bu_{i}^{'}(t_{n+1})+\bu_{i,h}^{'n}\|_{L^\infty}\|\nabla \bu_{j}(t_{n+1})\|\|\nabla \bphi_{j,h}^{n+1}\|\\&\le \frac{\talpha_j}{22}\|\nabla\bphi_{j,h}^{n+1}\|^2+\Delta t^6\frac{C}{\talpha_j}.
\end{align*}

Assume $\mu>1$, apply Cauchy-Schwarz, Young’s, and H\"older's  inequalities, definition \eqref{eddy-viscosity-new}, triangle inequality and Lemma \ref{lemma1},
	to have
    \begin{align}
		2\mu\Delta t\lp (\tl_h^{n})^2\nabla \bfeta_{j}^{n+1},\nabla\bphi_{j,h}^{n+1}\rp&= 2\mu\Delta t\left(\tl_h^{n}\nabla \bfeta_{j}^{n+1},\tl_h^{n}\nabla\bphi_{j,h}^{n+1}\right)\le 2\mu\Delta t\|\tl_h^{n}\nabla \bfeta_{j}^{n+1}\|\|\tl_h^{n}\nabla\bphi_{j,h}^{n+1}\|\nonumber\\&\le (\mu-1)\Delta t\|\tl_h^{n}\nabla\bphi_{j,h}^{n+1}\|^2+\frac{\mu^2}{\mu-1}\Delta t\|\tl_h^{n}\nabla \bfeta_{j}^{n+1}\|^2\nonumber\\&\le (\mu-1)\Delta t\|\tl_h^{n}\nabla\bphi_{j,h}^{n+1}\|^2+\frac{\mu^2}{\mu-1}\Delta t\|(\tl_h^{n})^2\|_\infty\|\nabla \bfeta_{j}^{n+1}\|^2
        \nonumber\\&\le (\mu-1)\Delta t\|\tl_h^{n}\nabla\bphi_{j,h}^{n+1}\|^2+\frac{\mu^2}{\mu-1}\Delta t\sum_{i=1}^J\|\bu_{i,h}^{'n}\|_\infty^2\|\nabla \bfeta_{j}^{n+1}\|^2\nonumber\\&\le (\mu-1)\Delta t\|\tl_h^{n}\nabla\bphi_{j,h}^{n+1}\|^2+\frac{C\Delta t}{\mu-1}\|\nabla \bfeta_{j}^{n+1}\|^2.
	\end{align}
    
Using Taylor’s series expansion, Cauchy-Schwarz inequality, non-linear bound in \eqref{nonlinearbound3}, regularity assumption, and Young's inequality, the last term of \eqref{orthdecomp} is bounded above as
\begin{align*}
    |-G(t,\bu_j,\bphi_{j,h}^{n+1})|\le\frac{\talpha_j}{22}\|\nabla\bphi_{j,h}^{n+1}\|^2+\Delta t^4\frac{C}{\talpha_j}\big(\|\bu_{j,ttt}(s_1^*)\|^2+\|\nabla\bu_{j,h}^{'n}\|^2\|\nabla\bu_{j,tt}(s_2^{*})\|^2\\+\|\bu_{j,tt}(s_3^{*})\|^2+\|\nu_j^{'}\|_\infty^2\|\nabla\bu_{j,tt}(s_4^*)\|^2\big),
\end{align*}
for some $s_i^*\in [t^{n-1},t^{n+1}]$, $i=\overline{1,4}$.

Using these estimates in \eqref{orthdecomp} and simplifying, provides
\begin{align}
    \frac{1}{4\Delta t}\Big(&\|\bphi_{j,h}^{n+1}\|^2
    +\|2\bphi_{j,h}^{n+1}-\bphi_{j,h}^{n}\|^2-\|\bphi_{j,h}^{n}\|^2-\|2\bphi_{j,h}^{n}-\bphi_{j,h}^{n-1}\|^2\Big)
    +\frac{\Bar{\nu}_{\min}}{2}\|\nabla\bphi_{j,h}^{n+1}\|^2
+\frac{\gamma}{2}\|\nabla\cdot\bphi_{j,h}^{n+1}\|^2\nonumber\\&+(\mu-1)\Delta t\|\tl_h^n\nabla \bphi_{j,h}^{n+1}\|^2+\left(\frac{1}{8\Delta t}-\frac{C}{\talpha_j}\|\nabla\cdot\bu_{j,h}^{'n}\|_{L^\infty}^2\right)\|\bphi_{j,h}^{n+1}-2\bphi_{j,h}^{n}+\bphi_{j,h}^{n-1}\|^2\nonumber\\
&\leq
    \|\nu_j^{'}\|_{\infty}\|\nabla 
		\bphi_{j,h}^{n}\|^2+\frac{1}{2}\|\nu_j^{'}\|_{\infty}\|\nabla 
		\bphi_{j,h}^{n-1}\|^2+\frac{44\|\nu_j^{'}\|_{\infty}^2}{\talpha_j}\|\nabla 
		\bfeta_{j}^{n}\|^2+\frac{11\|\nu_j^{'}\|_{\infty}}{\talpha_j}\|\nabla 
		\bfeta_{j}^{n-1}\|^2\nonumber\\&+\frac{C\mu^2\Delta t^2}{\talpha_j}\sum_{i=1}^J\left(\|\bphi_{i,h}^n\|^2+\|\bphi_{i,h}^{n-1}\|^2\right)+\frac{C\mu^2\Delta t^2}{\talpha_j}\sum_{i=1}^J\left(\|\bfeta_i^n\|^2+\|\bfeta_{i}^{n-1}\|^2\right)+\frac{C\Delta t}{\mu-1}\|\nabla \bfeta_{j}^{n+1}\|^2\nonumber\\&+\left(\frac{11\|\Bar{\nu}\|_\infty^2}{2\talpha_j}+C\gamma\right)\|\nabla\bfeta_{j}^{n+1}\|^2+\frac{C}{\talpha_j}\Big(\|\bphi_{j,h}^n\|^2+\|\bphi_{j,h}^{n-1}\|^2\Big)
  +\frac{C}{\talpha_j}\|\nabla(2\bu_{j,h}^n-\bu_{j,h}^{n-1})\|^2\|\nabla\bfeta_{j}^{n+1}\|^2\nonumber\\&+\frac{C}{\talpha_j}\|\nabla \bu_{j,h}^{'n}\|^2\|\nabla(\bfeta_{j}^{n+1}-2\bfeta_{j}^{n}+\bfeta_{j}^{n-1})\|^2+\frac{C}{\talpha_j}\|\nabla(2\bfeta_{j}^n-\bfeta_{j}^{n-1})\|^2+\Delta t^6\frac{C}{\talpha_j}\nonumber\\&+\Delta t^4\frac{C}{\talpha_j}\left(\|\bu_{j,ttt}(s_1^*)\|^2+\|\nabla\bu_{j,h}^{'n}\|^2\|\nabla\bu_{j,tt}(s_2^{*})\|^2+\|\bu_{j,tt}(s_3^{*})\|^2+\|\nu_j^{'}\|_\infty^2\|\nabla\bu_{j,tt}(s_4^*)\|^2\right).
\end{align}

Choose $\mu> 1$, time-step size given in \eqref{time-step-size}, drop the non-negative terms from the left-hand-side, triangle inequality, \eqref{alg-ineq} and rearrange, yields
\begin{align}
    &\frac{1}{4\Delta t}\Big(\|\bphi_{j,h}^{n+1}\|^2
    +\|2\bphi_{j,h}^{n+1}-\bphi_{j,h}^{n}\|^2-\|\bphi_{j,h}^{n}\|^2-\|2\bphi_{j,h}^{n}-\bphi_{j,h}^{n-1}\|^2\Big)
    +\frac{\Bar{\nu}_{\min}}{2}\left(\|\nabla\bphi_{j,h}^{n+1}\|^2-\|\nabla\bphi_{j,h}^{n}\|^2\right)\nonumber\\&+\left(\frac{\Bar{\nu}_{\min}}{2}-\|\nu_j^{'}\|_{\infty}\right)\left(\|\nabla\bphi_{j,h}^{n}\|^2-\|\nabla\bphi_{j,h}^{n-1}\|^2\right)+\frac{\talpha_j}{2}\|\nabla\bphi_{j,h}^{n-1}\|^2
\nonumber\\&\leq\frac{C}{\talpha_j}\left(\|\nabla 
		\bfeta_{j}^{n}\|^2+\|\nabla 
		\bfeta_{j}^{n-1}\|^2\right)+\frac{C\mu^2\Delta t^2}{\talpha_j}\sum_{i=1}^J\left(\|\bphi_{i,h}^n\|^2+\|\bphi_{i,h}^{n-1}\|^2+\|\bfeta_i^n\|^2+\|\bfeta_{i}^{n-1}\|^2\right)\nonumber\\&+\frac{C}{\talpha_j}\|\nabla \bu_{j,h}^{'n}\|^2\|\nabla(\bfeta_{j}^{n+1}-2\bfeta_{j}^{n}+\bfeta_{j}^{n-1})\|^2+\left(\frac{5\|\Bar{\nu}\|_\infty^2}{\talpha_j}+C\gamma\right)\|\nabla\bfeta_{j}^{n+1}\|^2\nonumber\\&+\frac{C}{\talpha_j}\Big(\|\bphi_{j,h}^n\|^2+\|\bphi_{j,h}^{n-1}\|^2\Big)
  +\frac{C}{\talpha_j}\|\nabla(2\bu_{j,h}^n-\bu_{j,h}^{n-1})\|^2\|\nabla\bfeta_{j}^{n+1}\|^2+\frac{C\Delta t}{\mu-1}\|\nabla \bfeta_{j}^{n+1}\|^2\nonumber\\&+\Delta t^4\frac{C}{\talpha_j}\bigg(\Delta t^2+\|\bu_{j,ttt}(s_1^*)\|^2+\|\nabla\bu_{j,h}^{'n}\|^2\|\nabla\bu_{j,tt}(s_2^{*})\|^2+\|\bu_{j,tt}(s_3^{*})\|^2+\|\nu_j^{'}\|_\infty^2\|\nabla\bu_{j,tt}(s_4^*)\|^2\bigg).
\end{align}

Now, multiplying both sides by $4\Delta t$ and taking sum for $n=1,\cdots, M-1$, assuming $\|\bphi^0_j\|=\|\bphi^1_j\|=\|\nabla\bphi^0_j\|=\|\nabla\bphi^1_j\|=0$, using $\Delta tM=T$, \eqref{alg-ineq}, stability estimate and regularity assumptions, we get\begin{align}
   \|\bphi_{j,h}^{M}\|^2+ \|2\bphi_{j,h}^{M}-\bphi_{j,h}^{M-1}\|^2+2\talpha_j\Delta t\sum_{n=2}^{M}\|\nabla\bphi_{j,h}^{n}\|^2\nonumber\\\leq \frac{C}{\talpha_j}\left\{h^{2k}+\Delta t^4+\Delta t\sum_{n=2}^{M-1}\left(\sum_{i=1}^J\|\bphi_{i,h}^n\|^2+\|\bphi_{j,h}^n\|^2\right)\right\}.
\end{align}

Summing over $j=1,\cdots\hspace{-0.35mm},J$, dropping non-negative term from left-hand-side, and grouping, to get
\begin{align}
\sum_{j=1}^J\|\bphi_{j,h}^{M}\|^2+2\talpha_{\min}\Delta t\sum_{n=2}^{M}\sum_{j=1}^J\|\nabla\bphi_{j,h}^{n}\|^2\le \frac{C}{\talpha_{\min}}(h^{2k}+\Delta t^4)+\frac{C\Delta t}{\talpha_{\min}}\sum_{n=2}^{M-1}\sum_{j=1}^J\|\bphi_{j,h}^n\|^2.
\end{align}

Applying the discrete Gr\"onwall Lemma \ref{dgl}, and using $\Delta t M=T$, we have
\begin{align}
\sum_{j=1}^J\|\bphi_{j,h}^{M}\|^2+2\talpha_{\min}\Delta t\sum_{n=2}^{M}\sum_{j=1}^J\|\nabla\bphi_{j,h}^{n}\|^2\le \frac{C}{\talpha_{\min}}(h^{2k}+\Delta t^4).
\end{align}

Finally, again use the triangle and Young's inequalities to complete the proof.
\end{proof}

\begin{remark}
    The use of the Lemma \ref{lemma1} could lead to a mild CFL-like time-step restriction.
\end{remark}

\begin{remark}
    Theorem \ref{stability-BE-EEV}/Theorem \ref{stability-bdf-2} shows that for a pointwise divergence-free stable element, the BE-EEV/BDF-2-EEV scheme is unconditionally stable with respect to the time-step size. 
\end{remark}

\begin{remark}
   If the assumptions of Theorem \ref{stability-BE-EEV}/Theorem \ref{stability-bdf-2} holds, then equation \eqref{stability-couple-alg}/equation \eqref{stability-bdf2-statement}  shows $\|\nabla\cdot\bu_{j,h}^{n}\|\rightarrow 0$ as  $\gamma\rightarrow\infty$. That is, even with the weakly divergence-free stable pair, the BE-EEV/BDF-2-EEV scheme is unconditionally stable with respect to the time-step size for large $\gamma$.
\end{remark}\vspace{-1mm}

\section{Numerical experiments}\label{numerical-experiment}
In this section, we present a series of numerical tests that verify the predicted convergence rates and show the performance of the schemes on some benchmark problems. In all experiments, we consider $\bx=(x_1,x_2)$, and $\bx=(x_1,x_2, x_3)$ for 2D, and 3D problems, respectively. Unless othewise stated, all the simulations are carried out in Deal.II, which is a C++ finite element library, with $(\mathbb{Q}_2,\mathbb{Q}_1)$ TH element, and the total number of realization $J=20$. The direct solver UMFPACK \cite{davis2004algorithm} is used to solve the linear system.\vspace{-1mm}

\subsection{Convergence rates verification}

In the first experiment, we verify the theoretically found convergence rates beginning with the following 2D analytical solution:
\[ {\bu}=\left(\begin{array}{c} \cos x_2+(1+e^t)\sin x_2 \\ \sin x_1+(1+e^t)\cos x_1 \end{array} \right)\;\;\text{and}\;\; \ p =\sin(x_1+x_2)(1+e^t),
\]
on domain $\cD=[0,1]^2$. Then, we introduce noise as $\bu_j=(1+k_j\epsilon)\bu$, and $p_j=(1+k_j\epsilon)p$, where $\epsilon$ is a perturbation parameter, $k_j:=(-1)^{j+1}4\lceil j/2\rceil/J$, $j=1,2,\cdots, J$, unless otherwise stated. The analytical solution is clearly divergence-free. We consider $\epsilon=10^{-3}$, and $10^{-2}$, this will introduce noise in the initial and boundary conditions, and the forcing functions. The forcing function $\bif_j$ is computed using the above synthetic data into \eqref{gov1}. We assume the viscosity $\nu$ is a continuous uniform random variable, and consider three random samples of size $J$ with $\bE[\nu]=10^{-3},10^{-4}$, and $10^{-5}$ where the samples are generated with 10$\%$ variation from the mean. The boundary conditions are set as $\bu_{j,h}|_{\partial\cD}=\bu_j$, and initial conditions $\bu_{j,h}^0=\bu_{j}(\bx,0)$ for both the scheme and the second initial condition of the BDF-2-EEV schemes are set as $\bu_{j,h}^1=\bu_{j}(\bx,\Delta t)$. We define the error as $<\hspace{-1mm}\be\hspace{-1mm}>:=<\hspace{-1mm}\bu\hspace{-1mm}>-<\hspace{-1mm}\bu_{h}\hspace{-1mm}>$ for the mean velocity. We consider structured quadrilateral meshes, and $\mu=1$.

 \subsubsection{Spatial convergence}

To observe spatial convergence, we keep the temporal error small enough, and thus, fix a very short simulation end time $T=0.001$. We consider $\epsilon=10^{-2}$, and $10^{-3}$. We successively refine the mesh width $h$ by a factor of 1/2, run the simulations, record the errors and convergence rates in Table \ref{spatial-con-ep-0_001}-\ref{spatial-con-ep-0_01} of the BE-EEV scheme, and in Table \ref{bdf2-ep-0-01}-\ref{bdf2-ep-0-001} for the BDF-2-EEV scheme. We observe the second order spatial convergence rates of both schemes for all three random samples of viscosity with the different $\epsilon$, which support our theoretical finding for the $(\mathbb{Q}_2,\mathbb{Q}_1)$ element. 

    \begin{table}[!ht]
    \centering
    \small\begin{tabular}{|l|l|l|l|l|l|l|}
    \hline
        \multicolumn{7}{|c|}{BE-EEV scheme: Spatial convergence (fixed $T = 0.001,~~\Delta t=1/8$)}\\ \hline
        $\epsilon = 10^{-3}$ & \multicolumn{2}{c|}{$\bE[\nu] = 10^{-3}$} & \multicolumn{2}{c|}{$\bE[\nu] =10^{-4}$} & \multicolumn{2}{c|}{$\bE[\nu] = 10^{-5}$} \\ \hline
        $h$ & $\|\hspace{-1mm}<\hspace{-1mm}\be\hspace{-1mm}>\hspace{-1mm}\|_{2,1}$ & rate & $\|\hspace{-1mm}<\hspace{-1mm}\be\hspace{-1mm}>\hspace{-1mm}\|_{2,1}$ & rate & $\|\hspace{-1mm}<\hspace{-1mm}\be\hspace{-1mm}>\hspace{-1mm}\|_{2,1}$ & rate \\ \hline
        1/2 & 4.2263e-4 & --- & 4.2263e-4 & --- & 4.2263e-4 & --- \\ \hline
         1/4 & 1.0780e-4 & 1.97 & 1.0780e-4 & 1.97 & 1.0780e-4 & 1.97 \\ \hline
         1/8 & 2.7078e-5 & 1.99 & 2.7078e-5 & 1.99 & 2.7078e-5 & 1.99 \\ \hline
          1/16 & 6.7775e-6 & 2.00 & 6.7775e-6 & 2.00 & 6.7776e-6 & 2.00 \\ \hline
          1/32 & 1.6950e-6 & 2.00 & 1.6950e-6 & 2.00 & 1.6952e-6 & 2.00 \\ \hline
    \end{tabular}\vspace{-2mm}
    \caption{Spatial errors and convergence rates of the BE-EEV scheme with $\epsilon = 10^{-3}$, and $\gamma =$ 2.99e+7, 2.90e+7, and 2.96e+7 for $\bE[\nu] =10^{-3},10^{-4}$, and $10^{-5}$, respectively.}\label{spatial-con-ep-0_001}
\end{table}

 \begin{table}[!ht]
    \centering
    \small\begin{tabular}{|l|l|l|l|l|l|l|}
    \hline
        \multicolumn{7}{|c|}{BE-EEV scheme: Spatial convergence (fixed $T = 0.001,~\Delta t=1/8$)}\\ \hline
        $\epsilon = 10^{-2}$ & \multicolumn{2}{c|}{$\bE[\nu] = 10^{-3}$} & \multicolumn{2}{c|}{$\bE[\nu] = 10^{-4}$ } & \multicolumn{2}{c|}{$\bE[\nu] = 10^{-5}$} \\ \hline
        $h$ & $\|\hspace{-1mm}<\hspace{-1mm}\be\hspace{-1mm}>\hspace{-1mm}\|_{2,1}$ & rate & $\|\hspace{-1mm}<\hspace{-1mm}\be\hspace{-1mm}>\hspace{-1mm}\|_{2,1}$ & rate & $\|\hspace{-1mm}<\hspace{-1mm}\be\hspace{-1mm}>\hspace{-1mm}\|_{2,1}$ & rate \\ \hline
       1/2 & 4.2263e-4 & --- & 4.2263e-4 & --- & 4.2263e-4 & --- \\ \hline
         1/4 & 1.0780e-4 & 1.97 & 1.0780e-4 & 1.97 & 1.0780e-4 & 1.97 \\ \hline
         1/8 & 2.7078e-5 & 1.99 & 2.7078e-5 & 1.99 & 2.7078e-5 & 1.99 \\ \hline
          1/16 & 6.7781e-6 & 2.00 & 6.7782e-6 & 2.00 & 6.7782e-6 & 2.00 \\ \hline
          1/32 & 1.7002e-6 & 2.00 & 1.7003e-6 & 2.00 & 1.7004e-6 & 2.00 \\ \hline
    \end{tabular}\vspace{-2mm}
\caption{Spatial errors and convergence rates of the BE-EEV scheme with $\epsilon = 10^{-2}$, and $\gamma =$ 3.10e+7, 2.95e+7, and 4.50e+7 for $\bE[\nu] =10^{-3},10^{-4}$, and $10^{-5}$, respectively.}\label{spatial-con-ep-0_01}
\end{table}
\vspace{-4ex}
\begin{table}[h!]
    \centering
    \small\begin{tabular}{|c|c|c|c|c|c|c|}
        \hline
         \multicolumn{7}{|c|}{BDF-2-EEV scheme: Spatial convergence (fixed $T = 0.001$, $\Delta t = T/8$)} \\
        \hline
        $\epsilon = 10^{-3}$
        & \multicolumn{2}{c|}{$\mathbb{E}[\nu] = 10^{-3}$} & \multicolumn{2}{c|}{$\mathbb{E}[\nu] = 10^{-4}$} & \multicolumn{2}{c|}{$\mathbb{E}[\nu] = 10^{-5}$} \\
        \hline
        $h$ & $\|\hspace{-1mm}<\hspace{-1mm}\be\hspace{-1mm}>\hspace{-1mm}\|{2,1}$ & rate & $\|\hspace{-1mm}<\hspace{-1mm}\be\hspace{-1mm}>\hspace{-1mm}\|{2,1}$ & rate & $\|\hspace{-1mm}<\hspace{-1mm}\be\hspace{-1mm}>\hspace{-1mm}\|_{2,1}$ & rate \\ \hline
        $1/2$ & 3.9419e-4 & --  & 3.9419e-4 & --  & 3.9419e-4 & --    \\
        \hline
        $1/4$  & 1.0076e-4 & 1.97& 1.0076e-4 & 1.97 & 1.0076e-4 & 1.97  \\
        \hline
        $1/8$  & 2.5326e-5 & 1.99& 2.5326e-5 & 1.99 & 2.5326e-5 & 1.99  \\
        \hline
        $1/16$  & 6.3397e-6 & 2.00& 6.3397e-6 & 2.00 & 6.3397e-6 & 2.00 \\
        \hline
        $1/32$  & 1.5855e-6 & 2.00& 1.5970e-6 & 1.99 & 1.5855e-6 & 2.00  \\
        \hline
    \end{tabular}\vspace{-2mm}
    \caption{Spatial errors and convergence rates of the BDF-2-EEV scheme with $\gamma = $ 7.50e+4.} \vspace{-4ex}\label{bdf2-ep-0-001}
\end{table}

\begin{table}[h!]
    \centering
    \small\begin{tabular}{|c|c|c|c|c|c|c|}
        \hline
         \multicolumn{7}{|c|}{BDF-2-EEV scheme: Spatial convergence (fixed $T = 0.001$, $\Delta t = T/8$)} \\
        \hline
        $\epsilon = 10^{-2}$
        & \multicolumn{2}{c|}{$\mathbb{E}[\nu] = 10^{-3}$} & \multicolumn{2}{c|}{$\mathbb{E}[\nu] = 10^{-4}$} & \multicolumn{2}{c|}{$\mathbb{E}[\nu] = 10^{-5}$} \\
        \hline
        $h$ & $\|\hspace{-1mm}<\hspace{-1mm}\be\hspace{-1mm}>\hspace{-1mm}\|{2,1}$ & rate & $\|\hspace{-1mm}<\hspace{-1mm}\be\hspace{-1mm}>\hspace{-1mm}\|{2,1}$ & rate & $\|\hspace{-1mm}<\hspace{-1mm}\be\hspace{-1mm}>\hspace{-1mm}\|_{2,1}$ & rate \\ \hline
        $1/2$ & 3.9419e-4 & --  & 3.9419e-4 & --  & 3.9419e-4 & --    \\
        \hline
        $1/4$  & 1.0076e-4 & 1.97& 1.0076e-4 & 1.97 & 1.0076e-4 & 1.97  \\
        \hline
        $1/8$  & 2.5326e-5 & 1.99& 2.5326e-5 & 1.99 & 2.5326e-5 & 1.99  \\
        \hline
        $1/16$  & 6.3397e-6 & 2.00& 6.3397e-6 & 2.00 & 6.3397e-6 & 2.00 \\
        \hline
        $1/32$  & 1.5855e-6 & 2.00& 1.5861e-6 & 2.00 & 1.5860e-6 & 2.00  \\
        \hline
    \end{tabular}\vspace{-2mm}
    \caption{Spatial errors and convergence rates of the BDF-2-EEV scheme with $\gamma = $ 1.00e+5.}\vspace{-7ex}\label{bdf2-ep-0-01}
\end{table}

\subsubsection{Temporal convergence} 
 On the other hand, to observe temporal convergence, we keep fixed mesh size $h=1/64$ and simulation end time $T=1$. We run the simulations with various time-step sizes $\Delta t$ beginning with $T/1$, and $T/2$ for the BE-EEV, and BDF-2-EEV scheme, respectively and successively reduce it by a factor of 1/2, record the errors, compute the convergence rates, and present them in Table \ref{temp-con-ep-0_3}-\ref{tab:table_4.1.1.1}. The outcomes of the BE-EEV scheme is presented in Table \ref{temp-con-ep-0_3} which shows the temporal convergence rates approximately equal to 1. Since the linearized backward-Euler formula is used in the proposed SPP-EEV scheme, the found temporal convergence rate is optimal and in excellent agreement with the theory for all three samples of the random viscosity. Also, the Table \ref{tab:table_4.1.1.1} represents the outcomes of the BDF-2-EEV scheme and we observe the second order convergence rate, which is consistent with the theoretical predication.\vspace{-3mm}

\begin{table}[!ht]
    \centering
    \small\begin{tabular}{|l|l|l|l|l|l|l|}
    \hline
        \multicolumn{7}{|c|}{BE-EEV scheme: Temporal convergence (fixed $T = 1,~h=1/64$)}\\ \hline
        $\epsilon = 10^{-3}$  & \multicolumn{2}{c|}{$\bE[\nu] = 10^{-3}$} & \multicolumn{2}{c|}{$\bE[\nu] = 10^{-4}$} & \multicolumn{2}{c|}{$\bE[\nu] = 10^{-5}$} \\ \hline
        $\Delta t$ & $\|\hspace{-1mm}<\hspace{-1mm}\be\hspace{-1mm}>\hspace{-1mm}\|_{2,1}$ & rate & $\|\hspace{-1mm}<\hspace{-1mm}\be\hspace{-1mm}>\hspace{-1mm}\|_{2,1}$ & rate & $\|\hspace{-1mm}<\hspace{-1mm}\be\hspace{-1mm}>\hspace{-1mm}\|_{2,1}$ & rate \\ \hline
        $T/1$ & 7.4508e-1 & ---& 2.0552e-0 & --- & 2.6743e-0 & ---   \\ \hline
        $T/2$ & 3.0177e-1 & 1.30& 7.5781e-1 & 1.44 & 9.3285e-1 & 1.52  \\ \hline
        $T/4$ & 1.3384e-1 & 1.17& 3.2538e-1 & 1.22 & 3.9582e-1 & 1.24  \\ \hline
        $T/8$ & 6.2869e-2 & 1.09& 1.5124e-1 & 1.11 & 1.8335e-1 & 1.11   \\ \hline
        $T/16$ & 3.0489e-2 & 1.04& 7.3041e-2 & 1.05 & 8.8342e-2 & 1.05  \\ \hline
        $T/32$ & 1.5036e-2 & 1.02& 3.5951e-2 & 1.02 & 4.3415e-2 & 1.02  \\ \hline
        $T/64$ & 7.4896e-3 & 1.01& 1.7864e-2 & 1.01 & 2.1556e-2 & 1.01  \\ \hline
    \end{tabular}\vspace{-2mm}
     \caption{Temporal errors and convergence rates of the BE-EEV scheme with $\gamma=$ 1.00e+5.}\vspace{-7ex}\label{temp-con-ep-0_3}
\end{table}

\begin{table}[h!]
    \centering
    \small\begin{tabular}{|c|c|c|c|c|c|c|}
        \hline
         \multicolumn{7}{|c|}{BDF-2-EEV scheme: Temporal convergence (fixed $T = 1$, $h = 1/64$)} \\
        \hline
        $\epsilon = 10^{-3}$
        & \multicolumn{2}{c|}{$\bE[\nu] = 10^{-3}$} & \multicolumn{2}{c|}{$\bE[\nu] = 10^{-4}$} & \multicolumn{2}{c|}{$\bE[\nu] = 10^{-5}$} \\
        \hline
        $\Delta t$ & $\|\hspace{-1mm}<\hspace{-1mm}\be\hspace{-1mm}>\hspace{-1mm}\|_{2,1}$ & rate & $\|\hspace{-1mm}<\hspace{-1mm}\be\hspace{-1mm}>\hspace{-1mm}\|_{2,1}$ & rate & $\|\hspace{-1mm}<\hspace{-1mm}\be\hspace{-1mm}>\hspace{-1mm}\|_{2,1}$ & rate \\ \hline
        $T/2$   & 2.9730e-1 & --   & 2.9730e-1 & --   & 2.9730e-1 & --   \\ \hline
        $T/4$   & 9.3898e-2 & 1.66 & 9.5442e-2 & 1.64 & 9.5934e-2 & 1.63 \\ \hline
        $T/8$   & 2.5574e-2 & 1.88 & 2.6182e-2 & 1.87 & 2.6389e-2 & 1.86 \\ \hline
        $T/16$  & 6.4457e-3 & 1.99 & 6.6435e-3 & 1.98 & 6.7169e-3 & 1.97 \\ \hline
        $T/32$  & 1.4163e-3 & 2.19 & 1.4825e-3 & 2.16 & 1.5081e-3 & 2.16 \\ \hline
    \end{tabular}\vspace{-2mm}
    \caption{Temporal errors and convergence rates of the BDF-2-EEV scheme with $\epsilon = 10^{-3}$, $\gamma =$ 1.00e+5, and $k_j$ is the sample of a uniformly distributed random variable $k\sim U(-1,1)$.}\vspace{-6ex}
    \label{tab:table_4.1.1.1}
\end{table}


\subsection{Channel flow over a unit square step} 
\begin{figure} [ht]
		\centering	
  \subfloat[]{\includegraphics[width=0.49\textwidth,height=0.17\textwidth]{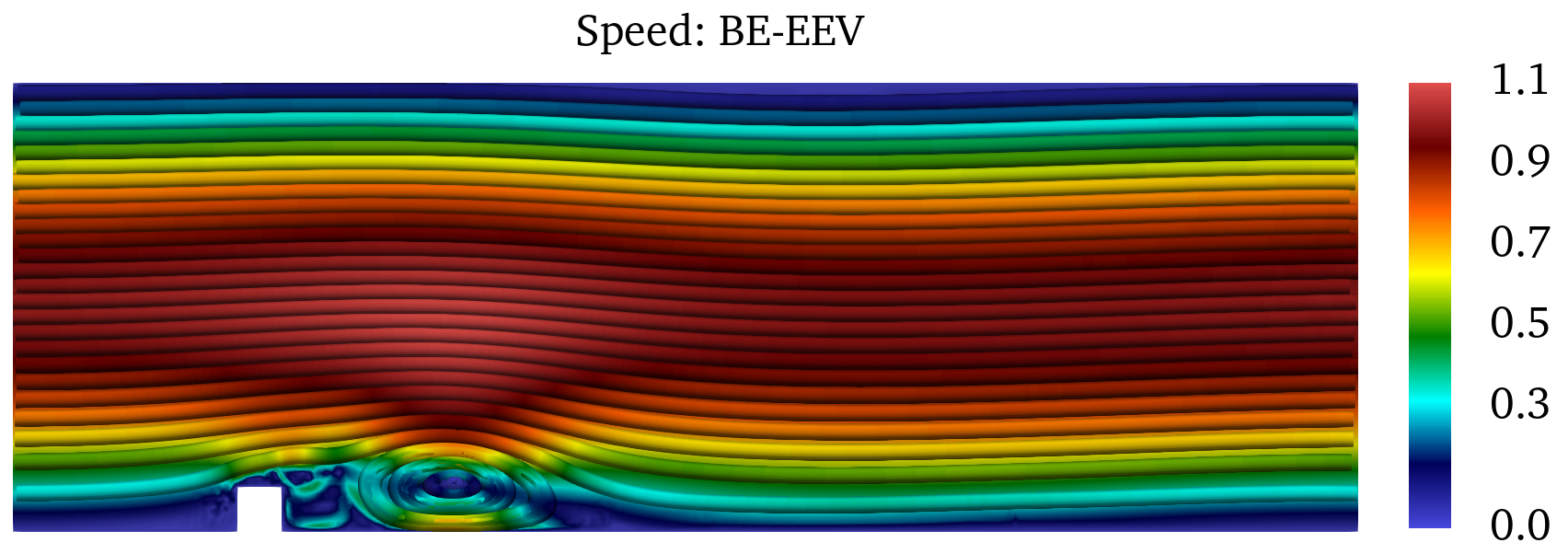}}\hspace{10mm}	
  \subfloat[]{\includegraphics[width=0.3\textwidth,height=0.17\textwidth]{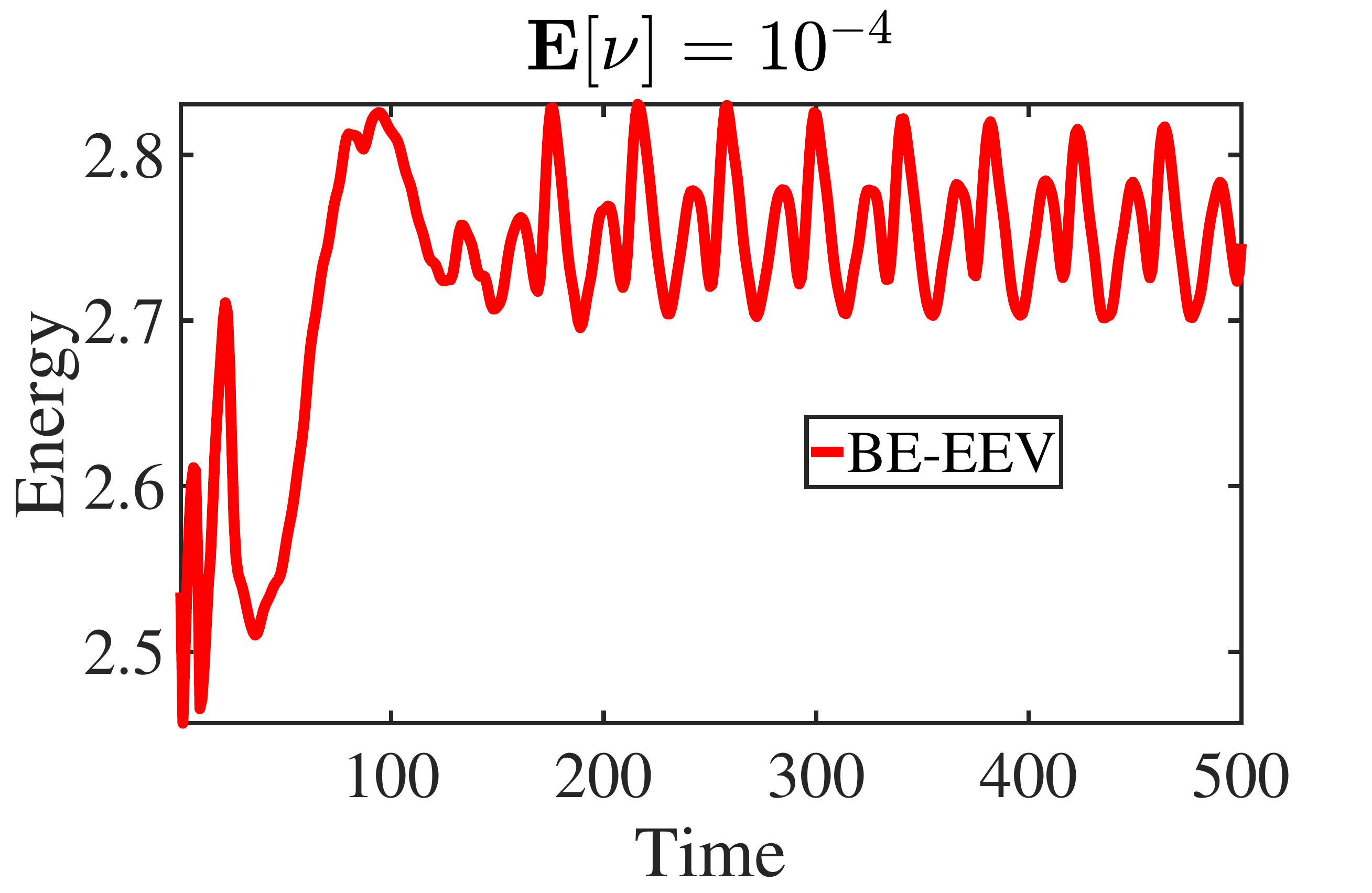}}\\\vspace{-3mm}
  \subfloat[]{\includegraphics[width=0.49\textwidth,height=0.17\textwidth]{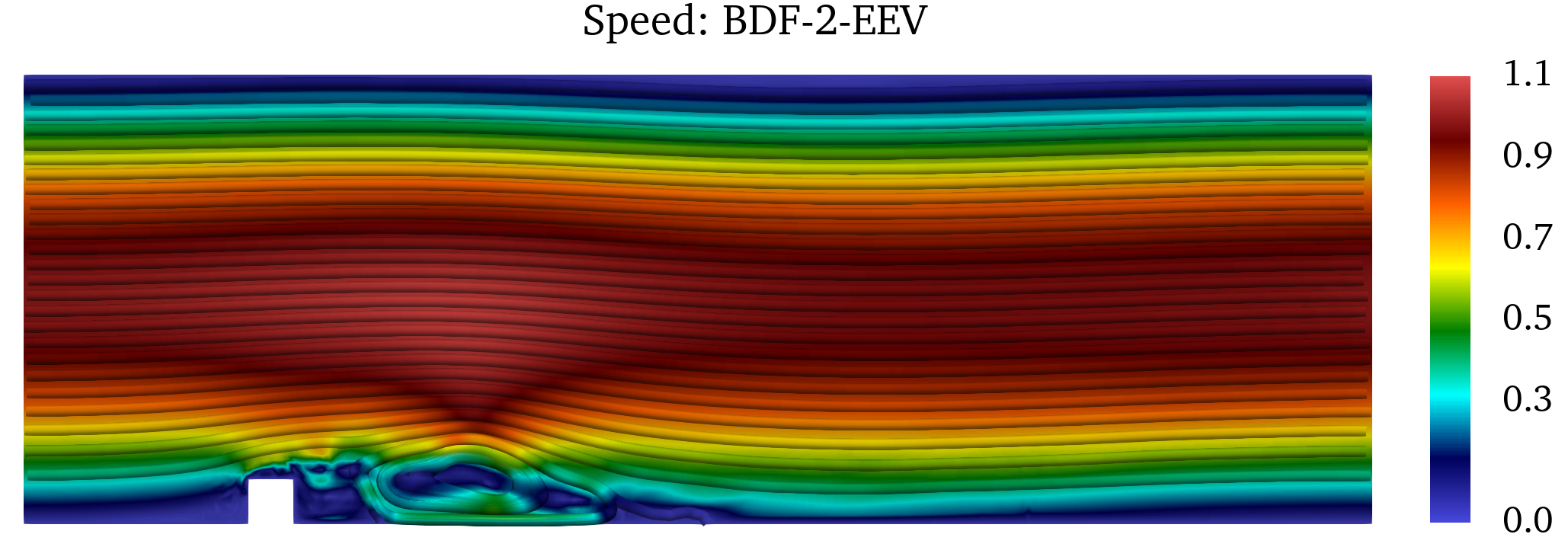}}\hspace{11mm}	
  \subfloat[]{\includegraphics[width=0.3\textwidth,height=0.2\textwidth]
  {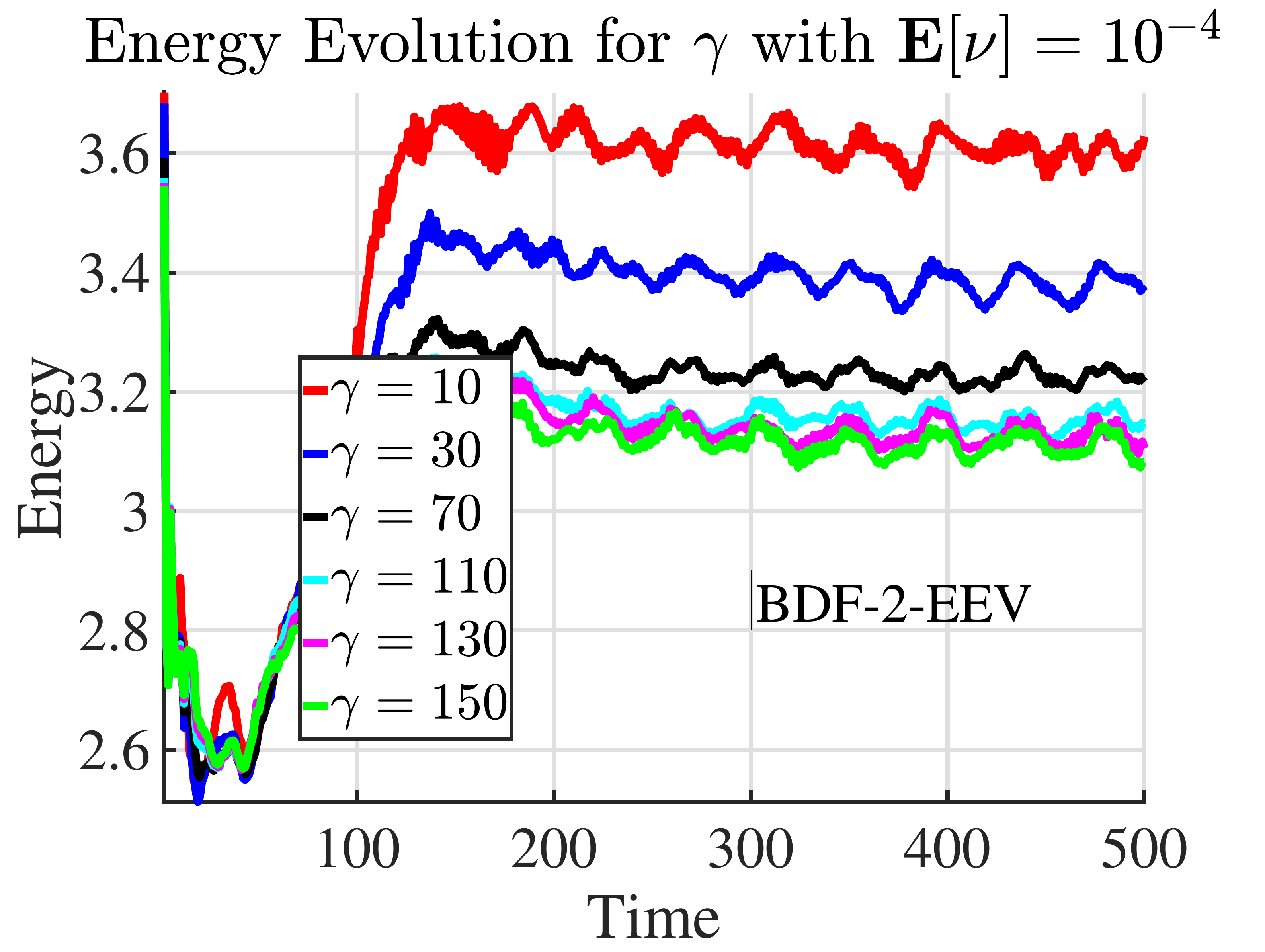}}\vspace{-6mm}
		\caption{\footnotesize{Flow over a step problem with $\mu=1$, $\epsilon=10^{-3}$, $\bE[\nu]=10^{-4}$, and dof $=100,835$: Streamlines over speed contour (of ensemble average solution) at $t=40$ with $\Delta t=0.08$ and $\gamma=10$ for the scheme (a) BE-EEV, and (c) BDF-2-EEV; plot of Energy vs. Time graph with $\Delta t=1$ for the scheme b) BE-EEV with $\gamma=10$, and (d) BDF-2-EEV with $\gamma$ varies.}}\vspace{-6ex}\label{channel-flow-over-step}
	\end{figure}

In this experiment, we consider a benchmark problem which is a 2D $30\times 10$ rectangular channel that has a $1\times 1$ step five units away from the inlet at the bottom. The following parabolic noisy flow is considered as the initial condition, inflow and outflow
\begin{align*}
    \bu_{j,h}(\bx,0)=\bu_{j,h}(\bx,t)|_{inlet,\; outlet}=\lp 1+k_j\epsilon\rp{{\frac{x_2(x_2-10)}{25}}\choose{0}},
\end{align*} where $k_j$ is the sample of a uniformly distributed random variable $k\sim U(-1,1)$, forcing $\bif_j=\textbf{0}$  for $j=1,2,\cdots,J$, and $\epsilon=10^{-3}$. To initialize the BDF-2-EEV scheme, we solve the BE-EEV scheme at the first time-step and use the solution as the second initial condition. No-slip boundary condition is applied to the domain walls and step. The channel geometry and unstructured quadrilateral meshes are generated using the free software Gmsh \cite{geuzaine2009gmsh}, which provides a total of 100,835 dof for each realization at each time-step. We run simulations in Deal.II using both the BE-EEV and BDF-2-EEV schemes with $\bE[\nu]=10^{-4}$ for a uniformly distributed random variable $\nu$. The outcomes are visualized in Paraview \cite{ayachit2015paraview} and are presented in Fig. \ref{channel-flow-over-step}. The speed contour plots in Fig. \ref{channel-flow-over-step} (a),(c), are at $t=40$ and are computed with the time-step size $\Delta t=0.08$; we observe stable solutions and these well agreement with the literature \cite{layton2008introduction}. On the other hand, the Energy vs. Time graphs in Fig. \ref{channel-flow-over-step} (b),(d), are computed for a time span of $[0,500]$ with $\Delta t=1$. The Fig. \ref{channel-flow-over-step} (b) represents the outcome of BE-EEV scheme with $\gamma=10$, and a long-range stable solution is observed. Similarly, in Fig. \ref{channel-flow-over-step} (d), a converging behavior of the energy plots is observed as $\gamma$ grows.\vspace{-2ex}

\begin{table}[h!]
    \centering
    \small\begin{tabular}{|c|c|c|c|c|}
    \hline
         \multicolumn{5}{|c|}{Fixed $T = 0.16$, $\Delta t=0.08$, $J=20$} \\
        \hline
        Scheme
        & Standard $O(\Delta t)$ & BE-EEV& Standard $O(\Delta t^2)$ & BDF-2-EEV \\
        \hline
        Time (s) & 156.38   & 67.29  & 157.64&  66.92\\ \hline
    \end{tabular}\vspace{-3mm}
    \caption{Comparison of computational time taken by standard and efficient schemes for two time-steps.}\vspace{-3ex}
    \label{comp:time}
\end{table}

To demonstrate the efficiency of the BE-EEV and BDF-2-EEV schemes over the standard schemes, we consider this channel flow over the unit step problem with the same mesh, $\Delta t=0.08$, end time $T=0.16$ and all other data remains the same. We compute the solving time for two time-steps of the efficient and traditional schemes and present them in Table \ref{comp:time}. The two standard schemes are equivalent to the first-order BE-EEV, and the second-order BDF-2-EEV schemes, but their system matrices are different for different realizations at each time-step. It is observed that the two efficient schemes require less computational time compared to the equivalent standard schemes, besides the saving of a huge memory requirement.

\subsection{2D RLDC problem and application of SCM}\label{RLDC}

Noises in the model parameters typically depend on high-dimensional random variables \cite{gunzburger2019evolve,gunzburger2014stochastic} and cause the model to suffer from the curse of dimensionality \cite{stoyanov2016dynamically}. Sampling-based methods such as MCM or multi-level MCM \cite{beck2020multilevel} require a significantly larger number of realizations compared to the recently introduced popular UQ approach, SCM \cite{babuvska2007stochastic, gunzburger2019evolve,gunzburger2014stochastic,stoyanov2013hierarchy,stoyanov2016dynamically}. Moreover, for realistic fluid flow simulation in engineering and geophysics, legacy codes are commonly used. SCM, was first introduced by Mathelin and Hussaini \cite{mathelin2003stochastic}, makes use of the regularity of the solution map concerning the random variables through the use of global polynomial approximation. Additionally, these schemes are independent of the PDE solvers, allowing them to be easily combined with any legacy codes. High-dimensional integration is performed by SCM using deterministic quadrature on sparse grids, whereas solving large-scale problems with high-dimensional random inputs requires the convergence of MCM, which generates unaffordable computational costs. To use SCM, it is common to assume $\nu_j(\bx)=\nu(\bx,\by_j(\omega)),$ and $\by_j(\omega)=(y_1(\omega),y_2(\omega),\cdots,y_N(\omega))\in\bGamma\subset\mathbb{R}^N$ be a finite $N\in\mathbb{N}$ dimensional vector \cite{babuvska2007stochastic,gunzburger2019evolve} distributed according to a joint probability density function $\rho(\by)$ in some parameter space $\bGamma=\prod\limits_{l=1}^N\Gamma_l$ with expectation $\bE[\by]=\textbf{0}$, and variance $Var[\by]=\textbf{I}_{N\times N}$ (the identity matrix) so that the statistical information about the QoI $\psi$ is approximated as $\bE[\psi(\bu)]=\int_{\bGamma}\psi(\bu,\by)\rho(\by)d\by\approx\sum\limits_{j=1}^Jw_j\psi(\bu,\by_j),$ for the quadrature weights $\{w_j\}_{j=1}^J$, and the number of realizations $J$ is the total number of stochastic collocation points. That is, the multi-dimensional integration is approximated by a one-dimensional summation only. We employ the Clenshaw–Curtis sparse grid \cite{stoyanov2016dynamically} as the SCM in this experiment, where the above quadrature weights and the high-dimensional stochastic collocation points are generated by Tasmanian \cite{stoyanov2015tasmanian}. The affinely independent viscosity $\nu_j(\bx)$ will be computed using the SCP in Karhunen-Lo\'eve expansion \cite{gunzburger2019evolve}.\begin{figure}[h!] 
	\begin{center}    
          \subfloat[]
		{\includegraphics[width = 0.3\textwidth, height=0.25\textwidth]{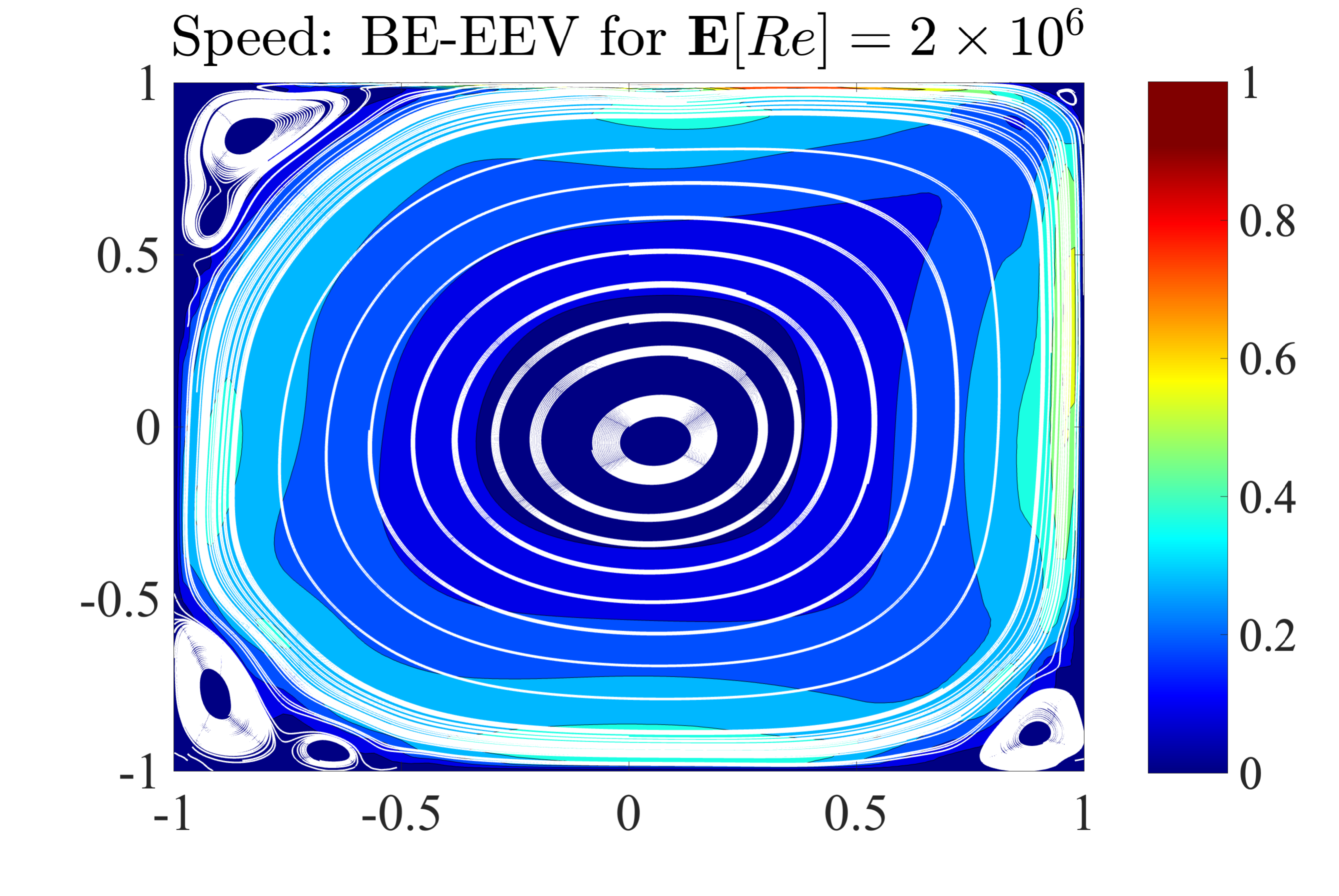}}\hspace{10mm} 
  \subfloat[]{
        \includegraphics[width = 0.35\textwidth, height=0.25\textwidth]{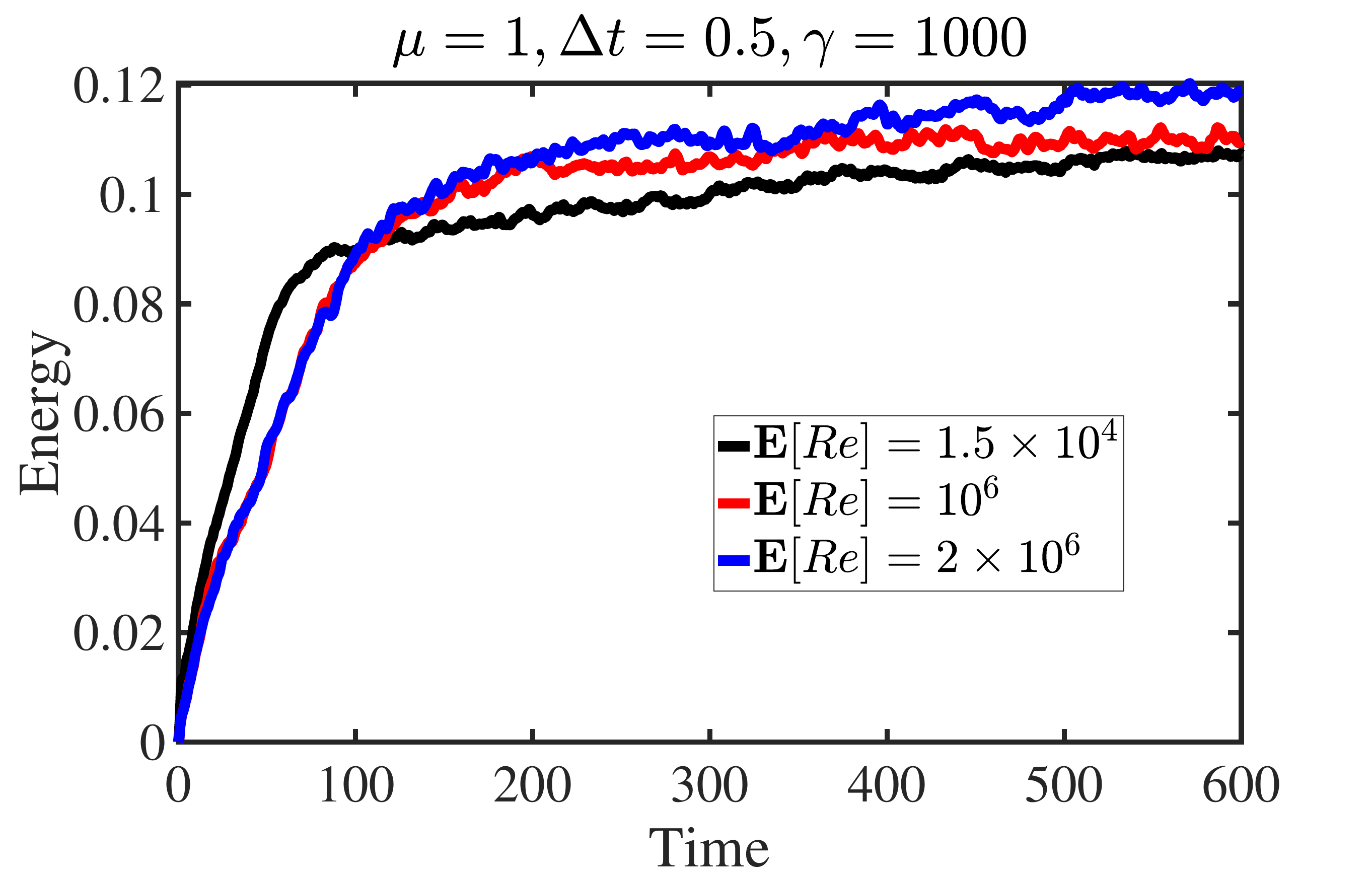}}
	\end{center}\vspace{-6mm}
	\caption{Variable 5D random viscosity in a RLDC problem with $\bE[Re]=2\times 10^6$: (a) Velocity solution (shown as streamlines over the speed contour) at $t=600$, and (b) Energy vs. Time plot for various $\bE[Re]$.}\vspace{-4ex}\label{RLDC_energy_curves}
\end{figure}

To this end, consider a 2D benchmark RLDC problem \cite{balajewicz2013low,fick2018stabilized,lee2019study, mohebujjaman2022efficient} with a domain $\cD=(-1,1)^2$. No-slip boundary conditions are applied to all sides except on the top wall (lid) of the cavity where we impose the following noise involved boundary condition:
\begin{align*}
\bu_{j,h}|_{lid}=\left (1+k_j\epsilon\right){{(1-x_1^2)^2}\choose{0}},
\end{align*}
so that the velocity of the boundary preserve the continuity. In this case, we model the sample viscosity using the Karhunen-Lo\'eve expansion:
\begin{align}
\nu_j=\nu({\bx}, {\by_j})=\frac{2}{\bE[Re]}\Bigg\{1+\left(\frac{\sqrt{\pi}l}{2}\right)^{\frac12}y_{j,1}(\omega)&+\sum_{i=1}^{2}\sqrt{\xi_i}\bigg(\sin\left(\frac{i\pi x_1}{2}\right)\sin\left(\frac{i\pi x_2}{2}\right)y_{j,2i}(\omega)\nonumber\\&+\cos\left(\frac{i\pi x_1}{2}\right)\cos\left(\frac{i\pi x_2}{2}\right)y_{j,2i+1}(\omega)\bigg)\Bigg\}, \label{eq:var-vis}
\end{align}
for a random vector ${\by_j}=(y_{j,1},y_{j,2},y_{j,3},y_{j,4},y_{j,5})\in\bGamma\subset\mathbb{R}^5$, where $\bGamma=[-\sqrt{3},\sqrt{3}]^5$, the correlation length $l=0.01$, and eigenvalues $$\sqrt{\xi_i}=(\sqrt{\pi}l)^{\frac12}exp\left(-\frac{(i\pi l)^2}{8}\right).$$
An unstructured bary-centered refined triangular mesh is considered that provides 364,920 dof per realization with $(\mathbb{P}_2,\mathbb{P}_1)$ element. We conducted the simulations on Freefem++ \cite{MR3043640} using the BDF-2-EEV scheme with the number of realizations $J=11$, end time $T=600$, step size $\Delta t=0.5$, $\gamma=1000$, $\mu=1$, $\epsilon=0.01$, $k\sim U(-1,1)$, and external force $\bif_j=\textbf{0}$. The computed solutions are then used in SCM to find the ensemble-weighted average of the quantity of interest. We represent the streamlines over the speed contour Fig. \ref{RLDC_energy_curves} (a) at $t=600$ for $\bE[Re]=2\times 10^6$ and the plot of Energy vs. Time graphs over the time period [0, 600] in Fig. \ref{RLDC_energy_curves} (b) for $\bE[Re]=1.5\times 10^4$, $10^6$, and $2\times 10^6$. The energy curves did not blow up and reached almost the statistical steady-state. We observe the formation of eddies in the corners of the cavity.\vspace{-2mm}

\subsection{3D RLDC problem} The RLDC problem has a domain $\cD=[-1,1]^3$. The flow is assumed to start moving from rest. The no-slip boundary conditions are applied to all surfaces except the lid, where we set the velocity \begin{align*}
    \bu_{j,h}|_{lid}=\lp 1+k_j\epsilon\rp <(2-x_1^2-x_2^2)^2/4,0,0>^T.
\end{align*}
\begin{figure}[h!] 
	\begin{center}    
          \subfloat[]
		{\includegraphics[width = 0.35\textwidth, height=0.35\textwidth]{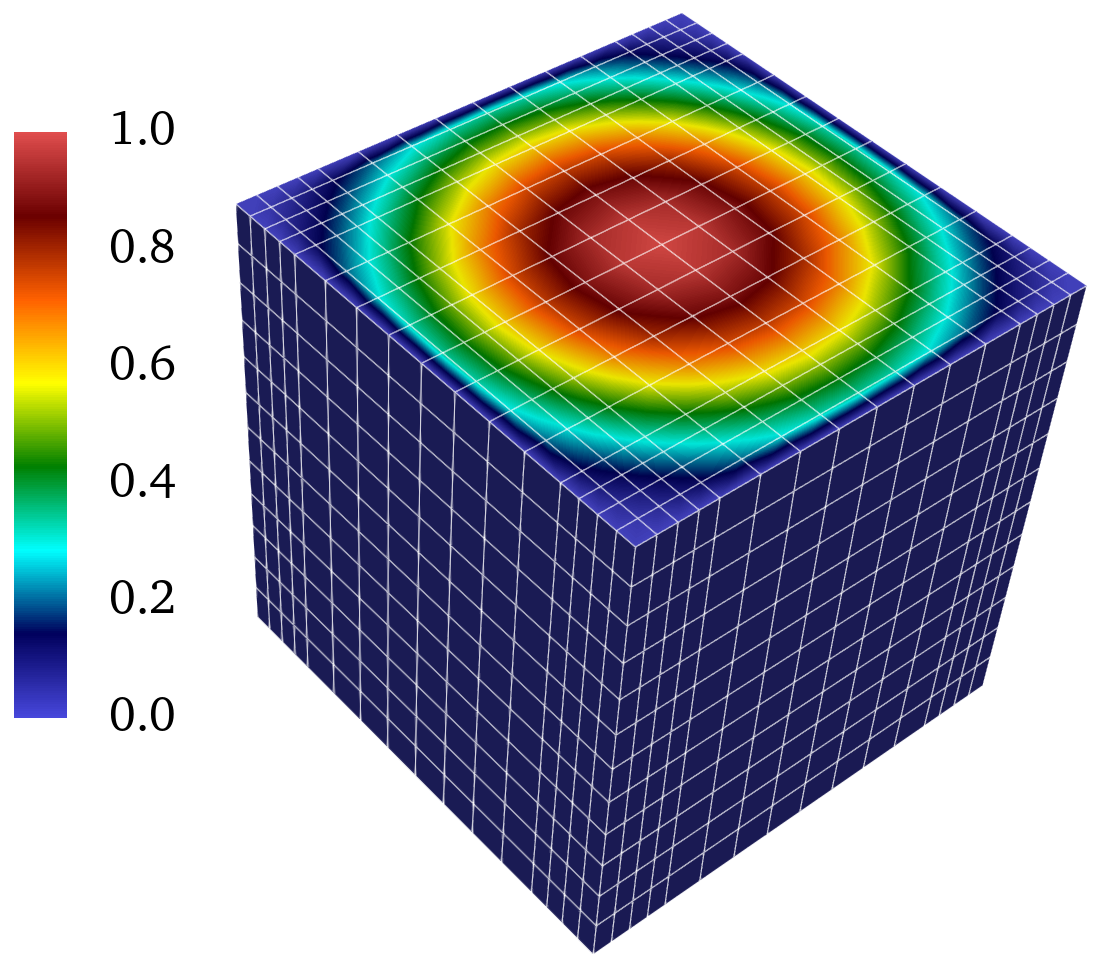}}\hspace{10mm} 
  \subfloat[]{
        \includegraphics[width = 0.35\textwidth, height=0.35\textwidth]{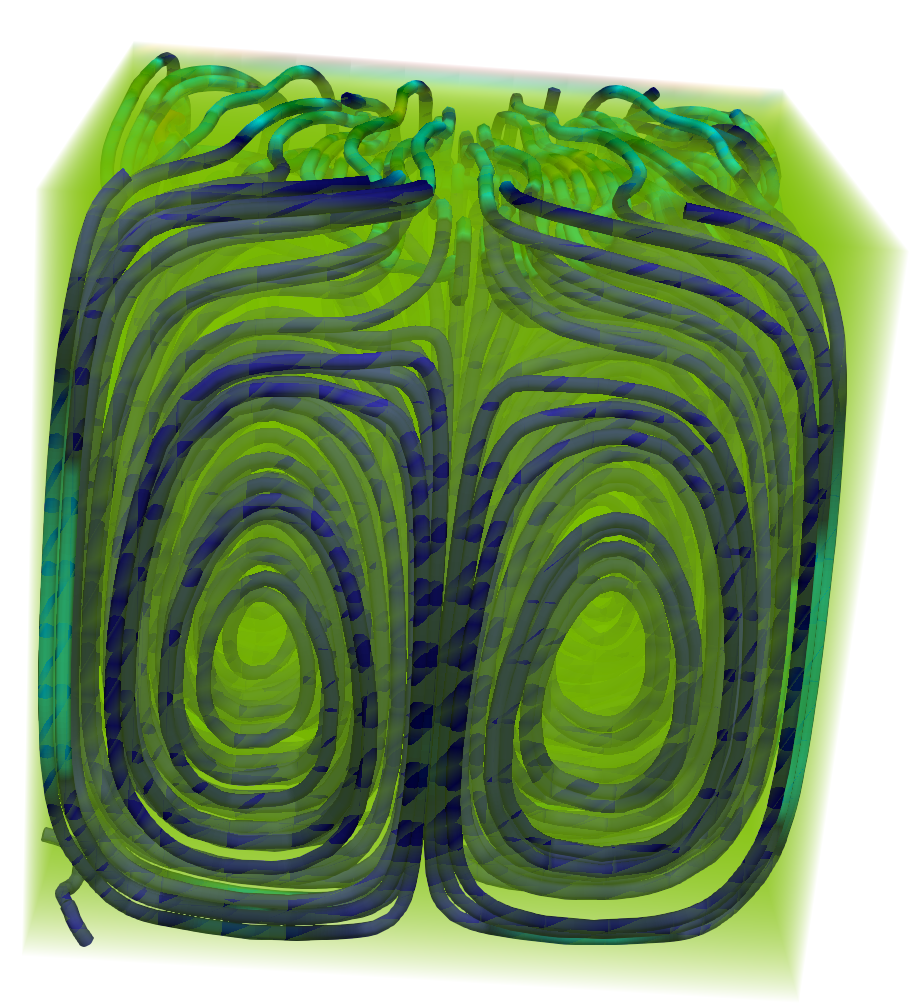}}
	\end{center}\vspace{-3mm}
	\caption{3D RLDC problem with $\mathbb{E}[Re]=20,000$  at $t=7$: (a) Top view of velocity solution (shown as speed contour over a sample computational mesh), and (b) Bottom view of streamlines over pressure contour.}\vspace{-5ex}\label{RLDC_3D}
\end{figure}
We assume the viscosity $\nu$ is a continuous uniform random variable, and consider a random sample of size $J=20$ so that $\bE[Re]=20,000$ where the sample is generated with 10$\%$ variation from the mean. The 3D geometry and its unstructured quadrilateral mesh are generated using Gmsh, which provides a total of 85432 dof per time-step for each realization. The following data are considered for this problem: $\Delta t=0.1$, $\mu=1$, $\gamma=100$, $(\mathbb{Q}_2,\mathbb{Q}_1)$ element, and $\epsilon=10^{-3}$. We run the simulations on Deal.II using the BDF-2-EEV scheme and visualize the solution at $t=7$ using Paraview. The speed contour over the computational mesh is shown in Fig. \ref{RLDC_3D} (a), and the bottom view of the streamlines over the pressure contour is shown in Fig. \ref{RLDC_3D}. The speed contour plot is clearly satisfying the imposed boundary condition, and streamlines shows symmetric pattern about the $xz$-plane.

\section{Conclusion and future research directions}\label{conclusion}

In this paper, we first have proposed a filter-based continuous EEV model for the turbulent SNSE problems. We then have proposed a novel family of efficient parameterized fully-discrete IMEX schemes for the EEV model assuming the presence of the random noises in the initial/boundary conditions, viscosity parameter, and forcing functions. The fully discrete schemes for the nonlinear SNSE are linearized; therefore, at each time-step, solving a nonlinear algebraic system is avoided, which dramatically reduces the number of arithmetic operations. Moreover, these schemes are designed in an elegant way so that at each time-step, the coefficient matrix remains the same for all realizations but different right-hand-side vectors. Thus, a single system matrix assembly is required per time-step and can take advantage of block linear solvers. Furthermore, a single LU decomposition or a single preconditioner can be built and reused for all realizations, allowing for significant savings in computational time and memory. The grad-div stabilization in the schemes enforces discrete divergence-free constraints pointwise, allowing the use of TH elements and larger time-step sizes.

We have analyzed, and tested two members of the proposed EEV family, namely, BE-EEV and BDF-2-EEV. The stability of the algorithms is proven rigorously and found them unconditionally stable with respect to the time-step size for large $\gamma$. The use of discrete inverse inequality in EEV models, which leads to suboptimal convergence for 3D problems, is avoided. It has been proven that optimal convergence, both in space and time, is achieved for both schemes. The temporal convergence of the BE-EEV and the BDF-2-EEV schemes is first- and second-order, respectively. The theoretical convergence rates are verified numerically using manufactured solutions for high expected Reynolds number problems. The schemes are implemented on benchmark problems: A 2D channel flow over a step problem, a 2D RLDC problem, and a 3D RLDC problem, and they have performed well. Both schemes perform well on the channel flow problem. In this case, the efficiency of these schemes is observed by comparing their computational time with that of the traditional equivalent schemes. Also, a converging trend of the BDF-2-EEV energy solutions is observed as $\gamma$ increases. The BDF-2-EEV turbulent scheme, in conjunction with the SCM, is implemented on the 2D RLDC problem and yields long-time stable solutions for $\bE[Re]=1.5\times 10^4$, $10^6$, and $2\times 10^6$. The BDF-2-EEV scheme is also tested on the 3D RLDC with $\bE[Re]=20,000$ and found a stable streamlines solution. 

The use of uniform boundedness of the discrete solutions could lead to a mild CFL-like time-step restriction; however, in the numerical experiments, it is not observed. These two EEV schemes could be enabling tools for large-scale multi-physics problems. A more efficient, second-order accurate, time-stepping penalty-projection algorithm will propose, analyze, and test for the turbulent flow problems as a future research work.

\section{Acknowledgment} The NSF is acknowledged for supporting this research through the grant DMS-2425308. We are grateful for the generous allocation of computing time provided by the Alabama Supercomputer Authority (ASA).
\bibliographystyle{plain}
\bibliography{BE}
\end{document}